\documentclass{article}
\usepackage[T2A]{fontenc}
\usepackage[cp1251]{inputenc}
\usepackage[english,russian]{babel}
\usepackage[tbtags]{amsmath}
\usepackage{amsfonts,amssymb,mathrsfs,amscd}
\numberwithin{equation}{section}
\newtheorem{remark}{Remark}[section]
\newtheorem{lemma}{Lemma}[section]
\newtheorem{theorem}{Theorem}[section]

\newtheorem{corollary}{Corollary}[section]

\renewcommand{\Re}{\mathop{\rm Re}}

\newcommand{\col}{\mathop{\rm col}}

\begin{document}
\begin{Large}
\thispagestyle{empty}
\begin{center}
{\bf Inverse spectral problem for a third-order differential operator\\
\vspace{5mm}
V. A. Zolotarev}\\

B. Verkin Institute for Low Temperature Physics and Engineering
of the National Academy of Sciences of Ukraine\\
47 Nauky Ave., Kharkiv, 61103, Ukraine

Department of Higher Mathematics and Informatics, V. N. Karazin Kharkov National University \\
4 Svobody Sq, Kharkov, 61077,  Ukraine

\end{center}
\vspace{5mm}

{\small {\bf Abstract.} Inverse spectral problem for a self-adjoint differential operator, which is the sum of the operator of the third derivative on a finite interval  and of the operator of multiplication by a real function (potential), is solved. Closed system of integral linear equations is obtained. Via solution to this system, the potential is calculated. It is shown that the main parameters of the obtained system of equations are expressed via spectral data of the initial operator. It is established that the potential is unambiguously defined by the four spectra.}
\vspace{5mm}


\vspace{5mm}

\begin{center}
{\bf Introduction}
\end{center}
\vspace{5mm}

Inverse problems were initiated by the work of V. A. Ambartsumian, 1929, and  were developed in works by G. Berg and N. Levinson \cite{1}. Seminal results in this field were obtained by V. A. Marchenko, I. M. Gelfand -- B. M. Levitan, L. D. Faddeyev \cite{1} -- \cite{3}. Inverse problems are divided into spectral (recovery of potential using spectral data) and scattering (recovery of potential using scattering data) problems. Use of transformation operators lies in the basis of the method for solving these problems for the Schrodinger operator. Efficient use of the method of inverse scattering problem made it possible to integrate several non-linear partial differential equations (Korteweg –- De Vries equation, etc.)

Generalization of these methods for a $p$th order differential operator ($p>2$) met with significant difficulties caused by the absence of transformation operators in this case. Construction of L - A pairs for non-linear equations describing oscillations in a dispersive medium \cite{4} -- \cite{8} (Camassa -- Holm, Degasperis -- Procesi equations, etc.) leads to a third-order operator $L$ (cubic string). Therefore, there arose necessity to solve inverse problems for differential operators of higher orders. Works \cite{9} -- \cite{12} solve inverse problems for the third-order operators. These works specify systems of functions that form the fundamental system of solutions (analogues of cosines and sines) for a third-order operator. These functions are the tool in terms of which inverse problems for the third-order differential operators are solved. Another method of the solution of inverse spectral problem ($p=3$) is given in \cite{13} -- \cite{15}.

This work solves inverse spectral problem for self-adjoint operators of the type
\begin{equation}
(L_q(\theta)y)(x)=iy'''(x)+q(x)y(x)\label{eq0.1}
\end{equation}
on a finite interval $x\in[0,l]$ ($0<l<\infty$) where $q(x)$ is a real function from $L^2(0,l)$ and domain of $L_q(\theta)$ is given by
\begin{equation}
\mathfrak{D}(L_q(\theta))=\{y\in W_2^3(0,l):y(0)=0,y'(l)=\theta y(0);y(l)=0\}\label{eq0.2}
\end{equation}
($\theta\in\mathbb{T}$). The work consists of 4 sections.

Section 1 studies spectral characteristics of the operator $L_0(\theta)$ ($q\equiv0$). It is shown that spectrum of such operator is simple, eigenfunctions of the operator $L_0(\theta)$ are constructed, and it is proved that they form an orthonormal basis in $L^2(0,l)$.

Section 2 studies the operator $L_q(\theta)$. A fundamental system of solution of this operator is constructed. This system is derived from the fundamental system of the operator $L_0(\theta)$ using analogues of transformation operators. Eigenfunctions of the operator $L_q(\theta)$ are found and resolvent of this operator is calculated.

Section 3 solves the direct problem, viz., the closed system of singular integral equations is obtained. Its solution allows to calculate potential. This system of equation is derived from the boundary value problem on the complex contour formed by the rays originating at zero with angle $2\pi/3$.

Section 4 solves the inverse problem. It is established that knowledge of four spectral data unambiguously allows to recover potential $q(x)$.

\section{Operator $L_0$}\label{s1}

{\bf 1.1} Equation
\begin{equation}
iD^3y(x)=\lambda^3y(x)\quad\left(d=\frac d{dx},\,\lambda\in\mathbb{C},\,x\in\mathbb{R}\right)\label{eq1.1}
\end{equation}
has three linearly independent solutions $\{e^{i\lambda\zeta_kx}\}_1^3$, where $\{\zeta_k\}_1^3$ are roots of the equation $\zeta^3=1$,
\begin{equation}
\zeta_1=1,\quad\zeta_2=\frac12(-1+i\sqrt3),\quad\zeta_3=\frac12(-1-i\sqrt3).\label{eq1.2}
\end{equation}
Any solution to equation \eqref{eq1.1} is a linear combination of these exponents. Instead of $\{e^{i\lambda\zeta_kx}\}_1^3$, it is convenient to chose another system of fundamental solutions to equation \eqref{eq1.1}, $\{s_p(i\lambda x)\}_0^2$, where
\begin{equation}
s_p(z)=\frac13\sum\limits_{k=1}^3\frac1{\zeta_k^p}e^{z\zeta_k}\quad(p=0,1,2)\label{eq1.3}
\end{equation}
which are analogues of cosines and sines for equation \eqref{eq1.1} \cite{9} -- \cite{11}.

The main properties of the functions $\{s_p(z)\}_0^2$ are listed in the following statement \cite{9} -- \cite{11}.

\begin{lemma}[\cite{9} -- \cite{11}]\label{l1.1}
Entire functions of exponential type $s_p(z)$ \eqref{eq1.3} have the following properties:

{\rm(i)} $s'_0(z)=s_2(z);\quad s'_1(z)=s_0(z);\quad s'_2(z)=s_1(z)\quad(y'(z)=dy(z)/dz);$

{\rm(ii)} $\overline{s_p(z)}=s_p(\overline{z})\quad(0\leq p\leq2)$;

{\rm(iii)} {\bf $p$-evenness}, $s_p(z\zeta_2)=\zeta_2^ps_p(z)\quad(0\leq p\leq2)$;

{\rm(iv)} Euler formula
$$e^{z\zeta_k}=s_0(z)+\zeta_ks_1(z)+\zeta_k^2s_2(z)\quad(1\leq k\leq3);$$

{\rm(v)} functions $\{s_p(z)\}$ \eqref{eq1.3} are solutions to the equation $y'''(z)=y(z)$ and satisfy the initial data,
$$s_0(0)=1,\quad s'_0(0)=0,\quad s''_0(0)=0;$$
$$s_1(0)=0,\quad s'_1(0)=1,\quad s''_1(0)=0;$$
$$s_2(0)=0,\quad s'_2(0)=0,\quad s''_2(0)=1;$$

{\rm(vi)} the main identity
$$s_0^3(z)+s_1^3(z)+s_2^3(z)-3s_0(z)s_1(z)s_2(z)=1;$$

{\rm(vii)} addition formulas
$$s_0(z+w)=s_0(z)s_0(w)+s_1(z)s_2(w)+s_2(z)s_1(w);$$
$$s_1(z+w)=s_0(z)s_1(w)+s_1(z)s_0(w)+s_2(z)s_2(w);$$
$$s_2(z+w)=s_0(z)s_2(w)+s_1(z)s_1(w)+s_2(z)s_0(w);$$

{\rm(viii)}
\begin{equation}
3s_0(z)s_0(w)=s_0(z+w)+s_0(z+\zeta_2w)+s_0(z+\zeta_3w);\label{eq1.4}
\end{equation}
$$3s_1(z)s_2(w)=s_0(z+w)+\zeta_2s_0(z+\zeta_2w)+\zeta_3s_0(z+\zeta_3w);$$
$$3s_1(z)s_0(w)=s_1(z+w)+s_1(z+\zeta_2w)+s_1(z+\zeta_3w);$$
$$3s_2(z)s_2(w)=s_1(z+w)+\zeta_2s_1(z+\zeta_2w)+\zeta_3s_1(z+\zeta_3w);$$
$$3s_2(z)s_0(w)=s_2(z+w)+s_2(z+\zeta_2w)+s_2(z+\zeta_3w);$$
$$3s_1(z)s_1(w)=s_2(z+w)+\zeta_3s_2(z+\zeta_2w)+\zeta_2s_2(z+\zeta_3w);$$

{\rm(ix)}
$$3s_0^2(z)=s_0(2z)+2s_0(-z);$$
$$3s_1^2(z)=s_2(2z)+2s_2(-z);$$
$$3s_2^2(z)=s_1(2z)+2s_1(-z);$$

{\rm(x)}
$$s_0^2(z)-s_1(z)s_2(z)=s_0(-z);$$
$$s_1^2(z)-s_0(z)s_2(z)=s_2(-z);$$
$$s_2^2(z)-s_0(z)s_1(z)=s_1(-z);$$

{\rm(xi)} Taylor formulas
$$s_0(z)=1+\frac{z^3}{3!}+\frac{z^6}{6!}+...;$$
$$s_1(z)=\frac z{1!}+\frac{z^4}{4!}+\frac{z^7}{7!}+...;$$
$$s_2(z)=\frac{z^2}{2!}+\frac{z^5}{5!}+\frac{z^8}{8!}+....$$
\end{lemma}

Solution to the Cauchy problem
\begin{equation}
iD^3y(x)=\lambda^3y(x)+f(x)\,(x\in\mathbb{R}_+);\quad y(0)=y_0,\quad y'(0)=y_1,\quad y''(0)=y_2\label{eq1.5}
\end{equation}
(for $f=0$), due to (v) \eqref{eq1.4}, is
\begin{equation}
y_0(\lambda,x)=y_0s_0(i\lambda x)+y_1\frac{s_1(i\lambda x)}{i\lambda}+y_2\frac{s_2(i\lambda x)}{(i\lambda)^2}.\label{eq1.6}
\end{equation}
Hence, using the arbitrary constant variation method, we obtain solution to the Cauchy problem \eqref{eq1.5} when $f\not=0$,
\begin{equation}
y(\lambda,x)=y_0(\lambda,x)-i\int\limits_0^x\frac{s_2(i\lambda(x-t))}{(i\lambda)^2}f(t)dt\label{eq1.7}
\end{equation}
where $y_0(\lambda,x)$ is given by \eqref{eq1.6}.

The straight lines
\begin{equation}
L_{\zeta_k}\stackrel{\rm def}{=}\{x\zeta_k:x\in\mathbb{R}\}\quad(1\leq k\leq3),\label{eq1.8}
\end{equation}
the direction vectors of which are the unit vectors $\{\zeta_k\}_1^3$, are equal to the union $L_{\zeta_k}=l_{\zeta_k}\cup\widehat{l}_{\zeta_k}$ where $l_{\zeta_k}$ are the half-lines (rays) originating at zero, and $\widehat{l}_{\zeta_k}$ are the half-lines ending at zero,
\begin{equation}
l_{\zeta_k}\stackrel{\rm def}{=}\{x\zeta_k:x\in\mathbb{R}_+\};\quad\widehat{l}_{\zeta_k}\stackrel{\rm def}{=}L_{\zeta_k}\backslash l_{\zeta_k}=\{x\zeta_k:x\in\mathbb{R}_-\}\quad(1\leq k\leq3).\label{eq1.9}
\end{equation}
These straight lines $\{L_{\zeta_k}\}_1^3$ \eqref{eq1.8} divide plane $\mathbb{C}$ into 6 sectors,
\begin{equation}
S_p\stackrel{\rm def}{=}\left\{z\in\mathbb{C}:\frac{2\pi}6(p-1)<\arg z<\frac{2\pi}6p\right\}\quad(1\leq p\leq 6).\label{eq1.10}
\end{equation}

\begin{lemma}[\cite{9}]\label{l1.2}
Zeros of the functions $\{s_p(z)\}_0^2$ \eqref{eq1.3} lie upon the rays $\{\widehat{l}_{\zeta_k}\}$ \eqref{eq1.9} and are equal to
\begin{equation}
\{-\zeta^lx_p(k)\}_{k=1}^\infty\label{eq1.11}
\end{equation}
where $l=-1$, $0$, $1$; and $x_p(k)$ are non-negative numbers enumerated in order of ascension of number $k$ ($p=0$, $1$, $2$). Numbers $x_0(k)$ (corresponding to $s_0(z)$) are positive roots of the equation
\begin{equation}
\cos\frac{\sqrt3}2x=-\frac12e^{-\frac32x}\quad(x_0(1)>0)\label{eq1.12}
\end{equation}
and $x_1(k)$, $x_2(k)$ (corresponding to $s_1(z)$, $s_2(z)$) are non-negative roots of the equations
\begin{equation}
\begin{array}{ccc}
{\displaystyle\cos\left(\frac{\sqrt3}2x-\frac\pi3\right)=\frac12e^{-\frac32x}}&(x_1(1)=0);\\
{\displaystyle\cos\left(\frac{\sqrt3}2x+\frac\pi3\right)=\frac12e^{-\frac32x}}&(x_2(1)=0).
\end{array}\label{eq1.13}
\end{equation}
The sequence $\{x_2(k)\}$ interlaces with the sequence $\{x_1(k)\}$ which, in its turn, interlaces with the sequence $\{x_0(k)\}$.
\end{lemma}

From equations \eqref{eq1.12}, \eqref{eq1.3}, it is easy to find the asymptotic behavior of $x_p(x)$ as $k\rightarrow\infty$ ($p=0$, 1, 2).

The rays $il_{\zeta_k}$ ($l_{\zeta_k}$ are given by \eqref{eq1.9}) divide plane $\mathbb{C}$ into the three sectors:
\begin{equation}
\begin{array}{ccc}
{\displaystyle\Omega_1\stackrel{\rm def}{=}\left\{\lambda\in\mathbb{C}:\frac{5\pi}6<\arg\lambda<\frac{11\pi}6\right\};\quad\Omega_2\stackrel{\rm def}{=}\left\{\lambda\in\mathbb{C}:\frac\pi2<\arg\lambda<\frac{5\pi}6\right\};}\\
{\displaystyle\Omega_3\stackrel{\rm def}{=}\left\{\lambda\in\mathbb{C}:-\frac\pi6<\arg\lambda<\frac\pi2\right\}.}
\end{array}\label{eq1.14}
\end{equation}
Obtain the asymptotic of functions $\{s_p(i\lambda x)\}_0^2$ in every sector $\{\Omega_k\}$ when $\lambda\rightarrow\infty$ ($x\in[0,l]$).

\begin{picture}(200,200)
\put(0,100){\vector(1,0){200}}
\put(100,0){\vector(0,1){200}}
\put(100,100){\vector(3,-1){100}}
\put(100,100){\vector(-3,-1){100}}
\qbezier(100,120)(130,110)(120,95)
\qbezier(119,94)(100,80)(81,94)
\qbezier(83,96)(80,120)(100,118)
\put(110,190){$il_{\zeta_2}$}
\put(97,170){$\circ$}
\put(110,170){$ix_p(k)$}
\put(130,120){$\Omega_3$}
\put(180,110){$l_{\zeta_1}$}
\put(180,62){$il_{\zeta_3}$}
\put(155,76){$\circ$}
\put(163,80){$i\zeta_3x_p(k)$}
\put(120,70){$\Omega_1$}
\put(6,77){$il_{\zeta_2}$}
\put(45,80){$\circ$}
\put(50,75){$i\zeta_2x_p(k)$}
\put(50,120){$\Omega_2$}
\end{picture}

\hspace{20mm} Fig. 1.

Equality
\begin{equation}
3s_0(i\lambda x)e^{-i\lambda\zeta_2x}=1+e^{i\lambda(\zeta_2-\zeta_1)x}+e^{i\lambda(\zeta_3-\zeta_1)x}\label{eq1.15}
\end{equation}
and identities
$$i\lambda(\zeta_2-\zeta_1)=\frac{\sqrt3}2[\beta\sqrt3-\alpha-i(\beta+\alpha\sqrt3)];$$
$$i\lambda(\zeta_3-\zeta_1)=\frac{\sqrt3}2[\beta\sqrt3+\alpha+i(\beta-\alpha\sqrt3)];$$
where $\lambda=\alpha+i\beta\in\mathbb{C}$ ($\alpha$, $\beta\in\mathbb{R}$), imply that, for $\beta\sqrt3-\alpha<0$ and $\beta\sqrt3+\alpha<0$, every exponent of the right-hand side of \eqref{eq1.15} is less than 1 in modulo and tends to zero when $\lambda\rightarrow\infty$ (inside the sector). So,
$$\lim\limits_{\lambda\rightarrow\infty}3s_0(i\lambda x)e^{-i\lambda\zeta_1x}=1\quad(\lambda\in\Omega_1).$$
Analogous considerations in the other sectors $\Omega_k$ \eqref{eq1.14} give the statement.

\begin{lemma}\label{l1.3}
For all $x\in[0,l]$, entire functions of exponential type $\{s_p(i\lambda x)\}_0^2$ have the following asymptotic in the sectors $\{\Omega_k\}_1^3$:
\begin{equation}
\lim\limits_{\lambda\rightarrow\infty}3\zeta_k^ps_p(i\lambda x)e^{-i\lambda\zeta_kx}=1\quad(\lambda\in\Omega_k)\label{eq1.16}
\end{equation}
($0\leq p\leq2$; $1\leq k\leq3$). Zeros of the functions $s_p(i\lambda x)$ are situated on the rays $\{il_{\zeta_k}\}_1^3$ separating the sectors $\{\Omega_k\}_1^3$.
\end{lemma}
\vspace{5mm}

{\bf 1.2.} In the space $L^2(0,l)$ ($0<l<\infty$), consider the self-adjoint operator $L_0(\theta)$ generated by the operation $iD^3$,
\begin{equation}
(L_0(\theta)y)(x)\stackrel{\rm def}{=}iD^3y(x)=iy'''(x),\label{eq1.17}
\end{equation}
the domain of which is
\begin{equation}
\mathfrak{D}(L_0(\theta))\stackrel{\rm def}{=}\{y\in W_2^3(0,l):y(0)=0;y'(l)=\theta y'(0);y(l)=0\}\quad(\theta\in\mathbb{T}).\label{eq1.18}
\end{equation}

\begin{remark}\label{r1.1}
The boundary conditions
$$y(0)=ih_0y''(0);\quad y'(l)=\theta y'(0);\quad y(l)=ih_ly''(l)\quad(h_0,h_l\in\mathbb{R},\,\theta\in\mathbb{T})$$
for the operation $iD^3$ also generate the self-adjoint operator $L_0(\theta,h_0,h_l)$ in the space $L^2(0,l)$ \cite{16,17}.
\end{remark}

The function
\begin{equation}
y_0(\lambda,x)=y_1\frac{s_1(i\lambda x)}{i\lambda}+y_2\frac{s_2(i\lambda x)}{(i\lambda)^2},\label{eq1.19}
\end{equation}
due to \eqref{eq1.6}, is a solution to equation \eqref{eq1.5} ($f=0$) and satisfies the boundary condition $y_0(\lambda,0)=0.$ The second and the third boundary conditions in
\eqref{eq1.18} give the system of linear equations
\begin{equation}
\left\{
\begin{array}{lll}
{\displaystyle y_1(s_0(i\lambda l)-\theta)+y_2\frac{s_1(i\lambda l)}{i\lambda}=0;}\\
{\displaystyle y_1\frac{s_1(i\lambda l)}{i\lambda}+y_2\frac{s_2(i\lambda l)}{(i\lambda)^2}=0.}
\end{array}\right.\label{eq1.20}
\end{equation}

\begin{remark}\label{r1.2}
For $\lambda=0$, system \eqref{eq1.20} becomes
$$\left\{
\begin{array}{lll}
y_1(1-\theta)+y_2l=0;\\
{\displaystyle y_1l+y_2\frac{l^2}2=0,}
\end{array}\right.$$
the determinant of which is ${\displaystyle\Delta=-\frac{l^2}2(\theta+1)}$. For $\theta\not=-1$, the determinant $\Delta\not=0$ and solution to the system is trivial, $y_1=y_2=0$, i. e., $y_0(0,x)\equiv0$. If $\theta=-1$, then the system has a non-zero solution and $y_0(0,x)=c(x^2-xl)$ ($c$ is a constant) is an eigenfunction of the operator $L_0(\theta)$ corresponding to the eigenvalue $\lambda^3=0$.
\end{remark}

System \eqref{eq1.20} has a non-trivial solution $y_1$, $y_2$ if and only if its determinant $\Delta_\theta(0,\lambda)$ vanishes, $\Delta_\theta(0,\lambda)=0$ where
$$\Delta_\theta(0,\lambda)\stackrel{\rm def}{=}-\frac1{(i\lambda)^2}[s_2(i\lambda l)(\theta-s_0(i\lambda l))+s_1^2(i\lambda l)]$$
and, due to (x) \eqref{eq1.4},
\begin{equation}
\Delta_\theta(0,\lambda)=-\frac1{(i\lambda)^2}[\theta s_2(i\lambda l)+s_2(-i\lambda l)].\label{eq1.21}
\end{equation}
Function $\Delta_\theta(0,\lambda)$ is said to be the {\bf characteristic function of the operator} $L_0(\theta)$ \eqref{eq1.17}, \eqref{eq1.18}. If $\lambda$ is a zero of $\Delta_\theta(0,\lambda)$ \eqref{eq1.21}, then $\lambda^3$ is an eigenvalue of the operator $L_0(\theta)$.

\begin{lemma}\label{l1.4}
The entire function of exponential type $\Delta(0,\lambda)$ \eqref{eq1.21} has the following properties:
\begin{equation}
\Delta_\theta(0,\lambda\zeta_2)=\Delta_\theta(0,\lambda);\quad\overline{\Delta_\theta(0,\lambda)}=\overline{\theta}\Delta_\theta(0,\overline{\lambda}).\label{eq1.22}
\end{equation}
Zeros of the function $\Delta_\theta(0,\lambda)$ are given by $\{\zeta_2^p\lambda_n(0,\theta)\}_{n=-\infty}^\infty$ ($p=0,$ $1$, $2$), where $\{\lambda_n(0,\theta)\}$ are simple real zeros of $\Delta_\theta(0,\lambda)$ enumerated in the ascending order, besides, ${\displaystyle\lambda_n(0,\theta)\in\left(\frac2l(\pi n+\varphi),\frac2l(\pi(n+1)+\varphi)\right)}$ ($n\in\mathbb{Z}$, $\theta=e^{2i\varphi}$, $\varphi\in[0,\pi]$) and the following asymptotic formulas are true:
\begin{equation}
\begin{array}{lll}
{\displaystyle\lambda_n(0,\theta)\rightarrow\frac2l\left(-\frac\pi6+\pi n+\varphi\right)+o\left(\frac1n\right)\quad(n\rightarrow\infty);}\\
{\displaystyle\lambda_n(0,\theta)\rightarrow\frac2l\left(\frac\pi6+\pi n+\varphi\right)+o\left(\frac1n\right)\quad(n\rightarrow-\infty).}
\end{array}\label{eq1.23}
\end{equation}
\end{lemma}

P r o o f. Relations \eqref{eq1.22} follow from (ii), (iii) \eqref{eq1.4}. According to \eqref{eq1.22}, it is sufficient to find location of the real zeros of $\Delta_\theta(0,\lambda)$. Hereinafter, we consider $\theta\not=\pm1$, i. e., $\lambda=0$ is not a zero of $\Delta_\theta(0,\lambda)$ (see Remark \ref{r1.2}).

Equation \eqref{eq1.21} implies that
\begin{equation}
\Delta_\theta(0,\lambda)=\frac{2e^{i\varphi}}{3\lambda^2}\sum\limits_k\zeta_k\cos(\lambda\zeta_kl+\varphi),\label{eq1.24}
\end{equation}
therefore equality $\Delta_\theta(0,\lambda)=0$ is equivalent to the relation
$$\cos(\lambda l+\varphi)-\cos\left(\frac{\lambda l}2-\varphi\right)\ch\frac{\lambda\sqrt3}2-\sqrt3\sin\left(\frac{\lambda l}2-\varphi\right)\sh\frac{\lambda l\sqrt3}2=0,$$
or
$$\cos\left(\frac{\lambda l}2-\varphi\right)\left(\cos\frac{3\lambda l}2-\ch\frac{\lambda l\sqrt3}2\right)=\sin\left(\frac{\lambda l}2-\varphi\right)\left(\sqrt3\sh\frac{\lambda\sqrt3l}2-\sin\frac{3\lambda l}2\right).$$
Since ${\displaystyle\ch\frac{\lambda\sqrt3l}2>\cos\frac{3\lambda l}2}$ ($\lambda\not=0$), the last equality can be rewritten as
\begin{equation}
\cot\left(\frac{\lambda l}2-\varphi\right)=f(\lambda);\quad f(\lambda)\stackrel{\rm def}{=}\frac{\displaystyle\sqrt3\sh\frac{\lambda\sqrt3l}2-\sin\frac{\displaystyle3\lambda l}2}{\displaystyle\cos\frac{3\lambda l}2-\ch\frac{\lambda\sqrt3l}2}\quad(\lambda\not=0).\label{eq1.25}
\end{equation}
Function $f(\lambda)$ is odd and $f(\lambda)>0$ for $\lambda\in\mathbb{R}_-$, $f(0)=0$, $f(\lambda)<0$ for $\lambda\in\mathbb{R}_+$. Note that
$$f'(\lambda)=l\cdot\frac{\displaystyle-3+3\ch\frac{\lambda l\sqrt3}2\cos\frac{3i\lambda}2+\sqrt3\sh\frac{\lambda l\sqrt3}2\sin\frac{3\lambda l}2}{\displaystyle\left(\cos\frac{3\lambda l}2-\ch\frac{\lambda l\sqrt3}2\right)^2}$$
and ${\displaystyle\cos\frac{3\lambda l}2<\ch\frac{\lambda l\sqrt3}2}$ ($\forall\lambda\in\mathbb{R}$, $\lambda\not=0$), therefore
$$-3+3\ch\frac{\lambda l\sqrt3}2\cos\frac{3\lambda l}2+\sqrt3\sh\frac{\lambda l\sqrt3}2\sin\frac{3\lambda l}2<-3+3\ch^2\frac{\lambda l\sqrt3}2+\sqrt3\sh\frac{\lambda l\sqrt 3}2\sin\frac{3\lambda l}2$$
$$=-\sqrt3\sh\frac{\lambda l\sqrt3}2\left(\sqrt3\sh\frac{\lambda l\sqrt3}2-\sin\frac{3\lambda l}2\right)<0\quad(\forall\lambda>0).$$
Hence, $f'(\lambda)<0$ for all $\lambda>0$, i. e., $f(\lambda)$ monotonically decreases on $\mathbb{R}_+$ from $0=f(0)$ to the asymptote $f(\lambda)\rightarrow-\sqrt3$ ($\lambda\rightarrow\infty$). On the semi-axis $\mathbb{R}_-$, the function $f(\lambda)$ also decreases (due to its oddness) and $f(\lambda)\rightarrow\sqrt3$ ($\lambda\rightarrow-\infty$).

Thus, the equation ${\displaystyle\cot\left(\frac{\lambda l}2-\varphi\right)=f(\lambda)}$ \eqref{eq1.25} in each of the intervals ${\displaystyle\left(\frac2l(\pi n+\varphi),\frac2l(\pi(n+1)+\varphi)\right)}$ has only one root $\lambda_n(0,\theta)$ ($n\in\mathbb{Z}$).

Asymptotic behavior \eqref{eq1.23} of the zeros $\lambda_n(0,\theta)$ when $n\rightarrow\pm\infty$ follows from \eqref{eq1.25} since $f(\lambda)\rightarrow\mp\sqrt3$ ($\lambda\rightarrow\pm\infty$), and thus ${\displaystyle\cot\left(\frac{\lambda l}2-\varphi\right)\rightarrow\mp\sqrt3}$ as $\lambda\rightarrow\pm\infty$.

The function $\Delta_\theta(0,\lambda)$ \eqref{eq1.21} cannot have complex zeros lying outside the straight lines $L_{\zeta_k}$ \eqref{eq1.8}. Because, if $w$ is such a zero, then $w^3$ is a complex eigenvalue of the self-adjoint operator $L_0(\theta)$, which is impossible. $\blacksquare$

To obtain the multiplicative expansion of $\Delta_\theta(0,\lambda)$, we use the Hadamard theorem \cite{18}.

\begin{theorem}[Hadamard \cite{18}]\label{t1.1}
Entire function of exponential type $f(z)$ can be written as
\begin{equation}
f(z)=ae^{bz}z^p\prod\limits_k\left(1-\frac z{z_k}\right)e^{\frac z{z_k}}\label{eq1.26}
\end{equation}
where $a$, $b\in\mathbb{C}$; $p\in\mathbb{Z}_+$; $\{z_k\}$ are zeros of $f(z)$ ($f(z_k)=0$, $z_k\not=0$, $\forall k$).
\end{theorem}

If $\theta\not=-1$ ($\lambda=0$ is not a zero of $\Delta_\theta(0,\lambda)$), then $p=0$ in the expansion \eqref{eq1.26} for $\Delta_\theta(0,\lambda)$. Every real root $\lambda_n(0,\theta)$ of the function $\Delta_\theta(0,\lambda)$ produces the series $\{\zeta_2^p\lambda_n(0,\theta)\}$ ($p=0$, 1, 2) (Lemma \ref{l1.3}), therefore
$$\left(1-\frac\lambda{\lambda_n(0)}\right)e^{\frac\lambda{\lambda_n(0)}}\cdot\left(1-\frac\lambda{\zeta_2\lambda_n(0)}\right)e^{\frac\lambda{\lambda_n(0)\zeta_2}}\cdot\left(1-
\frac\lambda{\zeta_3\lambda_n(0)}\right)e^{\frac\lambda{\lambda_n(0)\zeta_3}}=1-\frac{\lambda^3}{\lambda_n^3(0)},$$
and thus
$$\Delta_\theta(0,\lambda)=ae^{b\lambda}\prod\limits_n\left(1-\frac{\lambda^3}{\lambda_n^3(0,\theta)}\right).$$

Relation (xi) \eqref{eq1.4} implies that ${\displaystyle\Delta_\theta(0,0)=-\frac{l^2}2(\theta+1)}$, consequently, ${\displaystyle a=-\frac{l^2}2(\theta+1)}$. $b=0$ since ${\displaystyle\left.\left(\frac d{d\lambda}\Delta_\theta(0,\lambda)\right)\right|_{\lambda=0}=0}$ ($\theta+1\not=0$).

\begin{lemma}\label{l1.5}
The function $\Delta_\theta(0,\lambda)$ \eqref{eq1.21} is expressed via its real zeros $\{\lambda_n(0,\theta)\}$ by the formula
\begin{equation}
\Delta_\theta(0,\lambda)=\left\{
\begin{array}{lll}
{\displaystyle-\frac{l^2}2\prod\limits_n\left(1-\frac{\lambda^3}{\lambda_n^3(0,\theta)}\right)}\quad(\theta\not=-1);\\
{\displaystyle-\frac{2il^3}{5!}\overline{\lambda^3}\prod\limits_n\left(1-\frac{\lambda^3}{\lambda_n^3(0,\theta)}\right)}\quad(\theta=1).
\end{array}\right.\label{eq1.27}
\end{equation}
\end{lemma}
\vspace{5mm}

{\bf 1.3} Proceed to the eigenfunctions of operator $L_0(\theta)$. The second equality in \eqref{eq1.20} yields
$$y_1=-y_2\frac{s_2(i\lambda l)}{(i\lambda)^2}\cdot\left(\frac{s_1(i\lambda l)}{i\lambda}\right)^{-1}.$$
Substitute this expression into \eqref{eq1.19}, then
$$y_0(\lambda,x)=y_2\left(\frac{s_1(i\lambda l)}{i\lambda}\right)^{-1}\frac1{(i\lambda)^3}\{s_2(i\lambda x)s_1(i\lambda l)-s_1(i\lambda x)s_2(i\lambda l)\}.$$

\begin{remark}\label{r1.3}
If $s_1(i\lambda l)/i\lambda=0$ ($\lambda\not=0$), then a non-trivial solution to system \eqref{eq1.17} exists if $(a)$ $s_0(i\lambda l)=\theta$; $(b)$ $s_2(i\lambda l)=0$; $(c)$ $s_0(i\lambda l)=\theta$ and $s_2(i\lambda l)=0$. Show that these three cases are impossible for $\lambda\not=0$.

{\bf Case $(a)$.} If $s_1(i\lambda l)=0$ and $s_0(i\lambda l)=\theta$, then {\rm (iv)}, {\rm (vi)} \eqref{eq1.4} imply
$$s_2^3(i\lambda l)=1-\theta^3,\quad\zeta_p(e^{i\lambda\zeta_pl}-\theta)=s_2(i\lambda l)\quad(p=1,2,3)$$
Upon multiplying the last three equalities and using the first relation, we obtain
$$\theta(e^{i\lambda\zeta_1l}+e^{i\lambda\zeta_2l}+e^{i\lambda\zeta_3l})=e^{-i\lambda\zeta_1l}+e^{-i\lambda\zeta_2l}+e^{-i\lambda\zeta_3l},$$
and thus $\theta s_0(i\lambda l)=s_0(-i\lambda l)$, and since $s_0(i\lambda l)=\theta$, then $\theta^2=s_0(-i\lambda l)=\overline{\theta}$ ($\lambda\in\mathbb{R}$), consequently, $\theta^3=1$, i. e., $s_2(i\lambda l)=0$, which is impossible, because the roots of $s_1(z)$ and $s_2(z)$ don't intersect ($\lambda\not=0$, Lemma \ref{l1.2}).

{\bf Case $(b)$.} Equations $s_1(i\lambda l)=0$ and $s_2(i\lambda l)=0$, due to {\rm (iv)} \eqref{eq1.4}, imply that $s_0(i\lambda l)=e^{i\lambda\zeta_pl}$ ($p=1$, $2$, $3$), which is possible only when $\lambda=0$.

{\bf Case $(c)$.} If $s_0(i\lambda l)=\theta$ and $s_1(i\lambda l)=0$, $s_2(i\lambda l)=0$, then again {\rm (iv)} \eqref{eq1.4} yields that $\theta=e^{i\lambda\zeta_pl}$ ($p=1$, $2$, $3$) which is true only if $\lambda=0$.
\end{remark}

\begin{theorem}\label{t1.2}\label{t1.2}
Spectrum $\sigma(L_0)$ of the operator $L_0(\theta)$ \eqref{eq1.17} -- \eqref{eq1.19} is simple and
\begin{equation}
\sigma(L_0)=\{\lambda_n^3(0,\theta):n\in\mathbb{Z}\}\label{eq1.28}
\end{equation}
where $\{\lambda_n(0,\theta)\}$ are real zeros of characteristic function $\Delta_\theta(0,\lambda)$ \eqref{eq1.27} (Lemma \ref{l1.5}.

Eigenfunctions of the operator $L_0(\theta)$ corresponding to $\lambda=\lambda_n(0,\theta)$ are
\begin{equation}
\psi_n(0,\lambda,x)\stackrel{\rm def}{=}\frac1{a_n(\lambda)}\left\{\frac{s_2(i\lambda x)}{(i\lambda)^2}\cdot\frac{s_1(i\lambda l)}{i\lambda}-\frac{s_1(i\lambda x)}{i\lambda}\cdot\frac{s_2(i\lambda l)}{(i\lambda)^2}\right\}\quad(\lambda=\lambda_n(0,\theta))\label{eq1.29}
\end{equation}
besides, $\psi_n(0,\lambda\zeta_2,x)=\psi_n(0,\lambda,x)$, and numbers $a_n(\lambda)$ are found from the condition $\|\psi_n(0,\lambda,x)\|_{L^2(0,l)}=1$.
\end{theorem}

When $\theta=-1$, the value $\lambda=0$ is a zero of $\Delta_\theta(0,\lambda)$ \eqref{eq1.21}, eigenfunction corresponding to this value is given in Remark \ref{r1.2}.

\begin{remark}\label{r1.4}
For all $\lambda$, $w\in\mathbb{C}$, the following identity holds:
\begin{equation}
s_2(z)s_1(\lambda)-s_1(z)s_2(\lambda)=\frac1{i\sqrt3}\{s_0(z+\zeta_2w)-s_0(z+\zeta_3w)\}.\label{eq1.30}
\end{equation}
\end{remark}

Really, (viii) \eqref{eq1.4} implies
$$3(s_2(z)s_1(\lambda)-s_1(z)s_2(\lambda))=\zeta_2(s_0(w+\zeta_2z)-s_0(z+\zeta_2w))$$
$$+\zeta_3(s_0(w+\zeta_3z)-s_0(z+\zeta_3w))=(\zeta_2-\zeta_3)[s_0(z+\zeta_3w)-s_0(z+\zeta_2w)],$$
due to (iii) \eqref{eq1.4}, this gives \eqref{eq1.30} ($\zeta_2-\zeta_3=i\sqrt3$).

Using \eqref{eq1.30}, we obtain the following representations for the functions $\psi_n(0,\lambda,x)$ \eqref{eq1.29}:
$$\psi_n(0,\lambda,x)=\frac1{\sqrt3a_n(\lambda)\lambda^3}\{s_0(i\lambda(x+\zeta_2l))-s_0(i\lambda(x+\zeta_3l))\}.$$
\vspace{5mm}

{\bf 1.4} Calculate the resolvent $R_{L_0}(\lambda^3)=(L_0(\theta)-\lambda^3I)^{-1}$ of the operator $L_0(\theta)$ \eqref{eq1.17}, \eqref{eq1.18} and let $y=R_{L_0}(\lambda^3)f$, then
\begin{equation}
iy'''(x)-\lambda^3y(x)=f(x)\label{eq1.31}
\end{equation}
where $y\in\mathfrak{D}(L_0)$ \eqref{eq1.18} and $f\in L^2(0,l)$. The function
\begin{equation}
y(\lambda,x)=y_1\frac{s_1(i\lambda x)}{i\lambda}+y_2\frac{s_2(i\lambda x)}{(i\lambda)^2}-i\int\limits_0^x\frac{s_2(i\lambda(x-t))}{(i\lambda)^2}f(t)dt\label{eq1.33}
\end{equation}
is the solution to equation \eqref{eq1.30} and $y(\lambda,0)=0$. The second and the third boundary conditions \eqref{eq1.18} give the system of linear equations
$$\left\{
\begin{array}{lll}
{\displaystyle y_1(s_0(i\lambda l)-\theta)+y_2\frac{s_1(i\lambda l)}{i\lambda}=i\int\limits_0^l\frac{s_1(i\lambda(l-t))}{i\lambda}f(t)dt;}\\
{\displaystyle y_1\frac{s_1(i\lambda l)}{i\lambda}+y_2\frac{s_2(i\lambda l)}{(i\lambda)^2}=i\int\limits_0^l\frac{s_2(i\lambda(l-t))}{(i\lambda)^2}f(t)dt}
\end{array}\right.$$
the solution to which is
$$y_1=\frac i{(i\lambda)^3\Delta_\theta(0,\lambda)}\int\limits_0^l[s_1(i\lambda(l-t))s_2(i\lambda l)-s_2(i\lambda(l-t))s_1(i\lambda l)]f(t)dt;$$
$$y_2=\frac i{(i\lambda)^2\Delta_\theta(0,\lambda)}\int\limits_0^l[s_2(i\lambda(l-t))(s_0(i\lambda l)-\theta)-s_1(i\lambda(l-t))s_1(i\lambda l)]f(t)dt.$$
Substituting these values of $y_1$, $y_2$ in \eqref{eq1.33}, we obtain
$$y(\lambda,x)=\frac i{(i\lambda)^4\Delta_\theta(0,\lambda)}\left\{\int\limits_0^l\{s_1(i\lambda x)[s_1(i\lambda(l-t))s_2(i\lambda l)-s_2(i\lambda(l-t))s_1(i\lambda l)]\right.$$
\begin{equation}
+s_2(i\lambda x)[s_2(i\lambda(l-t))(s_0(i\lambda l)-\theta)-s_1(i\lambda(l-t))s_1(i\lambda l)]\}f(t)dt\label{eq1.34}
\end{equation}
$$\left.+\int\limits_0^xs_2(i\lambda(x-t))[\theta s_2(i\lambda l)+s_2(-i\lambda l)]f(t)dt\right\}.$$
Using (vii) \eqref{eq1.4}, we simplify the expression
$$s_1(l-t)s_2(l)-s_2(l-t)s_1(l)=s_2(l)[s_0(l)s_1(-t)+s_0(l)s_0(-t)+s_2(l)s_2(-t)]$$
$$-s_1(l)[s_0(l)s_2(-t)+s_2(l)s_0(-t)+s_1(l)s_1(-t)]=s_1(-t)[s_0(l)s_2(l)-s_1^2(l)]$$
$$+s_2(-t)[s_2^2(l)-s_0(l)s_1(l)]=s_2(-t)s_0(-l)-s_1(-t)s_2(-l)$$
due to (x) \eqref{eq1.4}. Analogously,
$$s_2(l-t)s_0(l)-s_1(l-t)s_1(t)=s_2(-t)s_0(-l)-s_0(-t)s_2(-l).$$
So,
$$s_1(x)[s_2(-t)s_0(-l)-s_1(-t)s_2(-l)]+s_2(x)[s_2(-t)s_0(-l)-s_0(-t)s_2(-l)]$$
$$=s_2(-t)s_2(x-l)-s_2(-l)s_2(x-t),$$
therefore formula \eqref{eq1.33} becomes
\begin{equation}
\begin{array}{ccc}
{\displaystyle y(\lambda,x)=\frac i{(i\lambda)^4\Delta_\theta(0,\lambda)}\left\{\int\limits_0^l[s_2(i\lambda(x-l))s_2(-i\lambda t)-s_2(i\lambda(x-t))s_2(-i\lambda l)\right.}\\
{\displaystyle\left.-\theta s_2(t\lambda x)s_2(i\lambda(z-t))]f(t)dt+\int\limits_0^xs_0(i\lambda(x-t))(\theta s_2(i\lambda l)+s_2(-i\lambda l))f(t)dt\right\}.}
\end{array}\label{eq1.35}
\end{equation}

\begin{theorem}\label{t1.3}
Resolvent $R_{L_0}(\lambda^3)=(L_0(t)-\lambda^3I)^{-1}$ of the operator $L_0(\theta)$ \eqref{eq1.17}, \eqref{eq1.18} is
\begin{equation}
\begin{array}{ccc}
{\displaystyle(R_{L_0}(\lambda^3)f)(x)=\frac i{(i\lambda)^4\Delta_\theta(0,\lambda)}\left\{\int\limits_0^l[s_2i\lambda(x-l)s_2(-i\lambda t)-\theta s_2(i\lambda x)\right.}\\
{\displaystyle\times s_2(i\lambda(l-t))]f(t)dt+\theta s_2(i\lambda l)\int\limits_0^xs_2(i\lambda(x-t))f(t)dt-s_2(-i\lambda l)\int\limits_x^ls_2(i\lambda(x-t))}\\
\left.\times f(t)dt\right\}
\end{array}\label{eq1.36}
\end{equation}
where $f\in L^2(0,l)$ and $\Delta_\theta(0,\lambda)$ is the characteristic function \eqref{eq1.27} of the operator $L_0(\theta)$.
\end{theorem}
\vspace{5mm}

{\bf 1.5} Resolvent of a self-adjoint operator $A$ with purely discrete spectrum acting in a Hilbert space $H$ is \cite{13}
$$R_A(\lambda)=(A-\lambda I)^{-1}=\sum\limits_n\frac{E_n}{\lambda_n-\lambda}$$
where $\lambda_n$ are eigenvalues of the operator $A$ and $E_n$ are orthogonal projections onto the proper subspaces corresponding to $\lambda_n$. Besides,
$$E_n=s-\lim\limits_{\lambda\rightarrow\lambda_n}(\lambda_n-\lambda)R_A(\lambda)$$
and ${\displaystyle\sum\limits_nE_n=I}$ which signifies completeness of proper subspaces $E_nH$ in the space $H$. We use this consideration to prove completeness of the eigenfunctions $\{\psi_n(0,\lambda,x)\}$ \eqref{eq1.29} in $L^2(0,l)$.

\begin{theorem}\label{t1.4}
The following relation is true:
\begin{equation}
(E_nf)(x)\stackrel{\rm def}{=}\lim\limits_{\lambda\rightarrow\lambda_n}(\lambda_n^3-\lambda^3)(R_{L_0}(\lambda^3)f)(x)=\int\limits_0^lf(t)\overline{\psi_n(0,\lambda,t)}dt\psi_n(0,\lambda,x)
\label{eq1.37}
\end{equation}
where $\psi_n(0,\lambda,x)$ is given by \eqref{eq1.29}, $\lambda_n=\lambda_n(0,t)$ is a real zero of the characteristic function $\Delta_\theta(0,\lambda)$ \eqref{eq1.21}, and for $a_n=a_n(\lambda)$ the following representation holds:
\begin{equation}
a_n^2=\left\{
\begin{array}{lll}
{\displaystyle\frac{s_2(i\lambda_nt)-s_2(-i\lambda_nl)}{2i}\frac{l^2}{\lambda_n^3}\prod\limits_{p\not=n}\left(1-\left(\frac{\lambda_n}{\lambda_p}\right)^3\right)\quad(
\theta\not=-1);}\\
{\displaystyle-s_2(i\lambda_nl)\frac{l^{\overline s}}{60\lambda_n^2}\cdot\prod\limits_{p\not=n}\left(1-\left(\frac{\lambda_n}{\lambda_p}\right)^3\right)\quad(\theta=-1)}
\end{array}
\right.
\end{equation}\label{eq1.38}
\end{theorem}

P r o o f. For $\theta\not=-1$, \eqref{eq1.35} and \eqref{eq1.27} imply
\begin{equation}
\begin{array}{ccc}
{\displaystyle(E_nf)(x)=\frac2{i\lambda_n(\theta+1)l^2\prod_n}\int\limits_0^l[s_2(i\lambda_n(x-l))s_2(-i\lambda_nt)-s_2(i\lambda_n(x-t))s_2(-i\lambda_nl)}\\
{\displaystyle-\theta s_2(i\lambda_nx)s_2(i\lambda_n(l-t))]f(t)dt\quad\left(\prod\limits_n\stackrel{\rm def}{=}\prod\limits_{p\not=n}\left(1-\left(\frac{\lambda_n}{\lambda_p}\right)^3\right)\right).}
\end{array}
\end{equation}
Consider the function
\begin{equation}
F(\lambda,x,t)\stackrel{\rm def}{=}s_2(i\lambda(x-l))s_2(-i\lambda t)-s_2(i\lambda(x-t))s_2(-i\lambda l)-\theta s_2(i\lambda x)s_2(i\lambda(l-t))\label{eq1.40}
\end{equation}
assuming that $\lambda$ is a zero of the function $\Delta(0,\lambda)$, i. e.,
\begin{equation}
\theta s_2(i\lambda l)+s_2(-i\lambda l)=0.\label{eq1.41}
\end{equation}
Note that
$$F(\lambda,x,0)=-s_2(i\lambda x)s_2(-i\lambda l)-\theta s_2(i\lambda x)s_2(i\lambda l)=0,$$
due to \eqref{eq1.41}. Upon differentiating with respect to $t$, we obtain
$$\partial_tF(\lambda,x,t)=(-i\lambda)[s_2(i\lambda(x-l))s_1(-i\lambda t)-s_1(i\lambda(x-t))s_2(-i\lambda l)-\theta s_2(i\lambda x)s_1(i\lambda(l-t))],$$
hence it follows that
$$\partial_tF(\lambda,x,0)=(-i\lambda)[-s_1(i\lambda x)s_2(-i\lambda l)-\theta s_2(i\lambda x)s_1(i\lambda l)]=-i\lambda\theta u(\lambda,x)$$
where, in accordance with \eqref{eq1.41},
\begin{equation}
u(\lambda,x)\stackrel{\rm def}{=}s_1(i\lambda x)s_2(i\lambda l)-s_2(i\lambda x)s_1(i\lambda l).\label{eq1.42}
\end{equation}
Repeated differentiation with respect to $t$ gives
$$\partial_t^2F(\lambda,x,t)=(-i\lambda)^2[s_2(i\lambda(x-l))s_0(-i\lambda t)-s_0(i\lambda(x-t))s_2(-i\lambda l)-\theta s_2(i\lambda x)s_0(i\lambda(l-t))],$$
and thus
$$\partial_t^2F(\lambda,x,0)=(-i\lambda)^2[s_2(i\lambda(x-l))-s_2(i\lambda x)s_2(-i\lambda l)-\theta s_2(i\lambda x)s_0(i\lambda l)]$$
$$=(-i\lambda)^2\{s_1(i\lambda x)s_1(-i\lambda l)-s_2(i\lambda x)[\theta s_0(i\lambda l)-s_0(-i\lambda l)]\},$$
due to (vii) \eqref{eq1.4}. For all such $\lambda$ that \eqref{eq1.41} holds, the following equality is true:
\begin{equation}
\frac{s_1(-i\lambda l)}{s_2(i\lambda l)}=\frac{\theta s_0(i\lambda l)-s_0(-i\lambda l)}{s_1(i\lambda l)}(\stackrel{\rm def}{=}R).\label{eq1.43}
\end{equation}
Really, taking into account \eqref{eq1.41},
$$s_1(-i\lambda l)s_1(i\lambda l)=s_2(i\lambda l)(\theta s_0(i\lambda l)-s_0(-i\lambda l))=-s_2(-i\lambda l)s_0(i\lambda l)-s_2(i\lambda l)s_0(-i\lambda l),$$
which coincides with (vii) \eqref{eq1.4}. Using \eqref{eq1.43}, \eqref{eq1.42}, we obtain that
$$\partial_t^2F(\lambda,x,0)=(-i\lambda)^2Ru(\lambda,x).$$
So, function $F(\lambda,x,t)$ \eqref{eq1.40}, for all $x\in[0,l]$, is the solution to the following Cauchy problem:
$$i\partial_t^3F(\lambda,x,t)=-\lambda^3F(\lambda,x,t);$$
$$\left.F(\lambda,x,t)\right|_{t=0}=0;\quad\left.\partial_tF(\lambda,x,t)\right|_{t=0}=(-i\lambda)\theta u(\lambda,x);$$
$$\left.\partial_t^2F(\lambda,x,t)\right|_{t=0}=(-i\lambda)^2Ru(\lambda,x).$$
Hence, according to \eqref{eq1.6}, we find that
$$F(\lambda,x,t)=\theta u(\lambda,x)s_1(-i\lambda t)+Ru(\lambda,x)s_2(-i\lambda t)=u(\lambda,x)[\theta s_1(-i\lambda t)+Rs_2(-i\lambda t)],$$
therefore, due to \eqref{eq1.41} -- \eqref{eq1.43},
$$F(\lambda,x,t)=\frac1{s_2(i\lambda l)}[-s_1(-i\lambda t)s_2(-i\lambda l)+s_2(-i\lambda t)s_1(-i\lambda l)]=-\frac1{s_2(i\lambda l)}\overline{u(\lambda,t)}u(\lambda,x).$$
Equation \eqref{eq1.29} implies that $u(\lambda,x)=-a_n(\lambda)(i\lambda_n)^3\psi_n(0,\lambda,x)$, and thus
$$F(\lambda,x,t)=-\frac{a_n^2(\lambda)\lambda_n^6}{s_2(i\lambda l)}\overline{\psi_n(0,\lambda,t)}\psi_n(0,\lambda,x).$$
Upon substituting this expression into \eqref{eq1.36}, we obtain
$$(E_nf)(x)=\frac{-2a_n^2(\lambda)\lambda_n^5}{s_2(i\lambda_nl)i(\theta+1)l^2\prod\limits_n}\int\limits_0^lf(t)\overline{\psi_n(0,\lambda,t)}dt\psi_n(0,\lambda,x).$$
And since $E_n$ is an orthogonal projection, then
$$a_n^2=\frac{s_2(i\lambda_nl)-s_2(-i\lambda_nl)}{2i\lambda_n^5}\cdot l^2\prod_n,$$
which proves \eqref{eq1.38} for $\theta\not=-1$. Equation \eqref{eq1.38} for $\theta=-1$ is proved analogously. $\blacksquare$

\begin{corollary}\label{c1.1}
The eigenfunctions $\{\psi_n(0,\lambda,x)\}$ \eqref{eq1.29} ($\lambda=\lambda_n(0,t)$ are zeros of the characteristic function $\Delta_\theta(0,\lambda)$ \eqref{eq1.21}) of the operator $L_0(\theta)$ \eqref{eq1.17}, \eqref{eq1.18} form the orthonormal basis in the space $L^2(0,l)$.
\end{corollary}

\section{Operator $L_q$}\label{s2}

{\bf 2.1} In the space $L^2(0,l)$ ($0<l<\infty$), consider the self-adjoint operator $L_q(\theta)$, which is a perturbation of the operator $L_0(\theta)$ \eqref{eq1.17}, \eqref{eq1.18},
\begin{equation}
(L_q(\theta)y)(x)\stackrel{\rm def}{=}iy'''(x)+q(x)y(x)\label{eq2.1}
\end{equation}
where $q(x)$ is a real function from $L^2(0,l)$, domain $\mathfrak{D}(L_q)$ of the operator $L_q(\theta)$ coincides with the domain $\mathfrak{D}(L_0)$ \eqref{eq1.18}, $\mathfrak{D}(L_q)=\mathfrak{D}(L_0)$.

Let $y(\lambda,x)$ be the solution to the equation
\begin{equation}
iy''(x)+q(x)y(x)=\lambda^3y(x)\quad(\lambda\in\mathbb{C})\label{eq2.2}
\end{equation}
satisfying the initial data
\begin{equation}
y(\lambda,0)=0,\quad y'(\lambda,0)=y_1,\quad y''(\lambda,0)=y_2\quad(y_1,y_2\in\mathbb{C}).\label{eq2.3}
\end{equation}
Then the difference
\begin{equation}
z(\lambda,x)\stackrel{\rm def}{=}y(\lambda,x)-y_0(\lambda,x)\label{eq2.4}
\end{equation}
($y_0(\lambda,x)$ is given by \eqref{eq1.19}) is a solution to the non-homogenous equation
$$iz'''(\lambda,x)=\lambda^3z(\lambda,x)-q(x)y(\lambda,x)$$
and satisfies the zero initial data,
$$z(\lambda,0)=0,\quad z'(\lambda,0)=0,\quad z''(\lambda,0)=0.$$
Hence, due to \eqref{eq1.7}, we find that
\begin{equation}
\begin{array}{ccc}
{\displaystyle z(\lambda,x)=i\int\limits_0^xK_1(\lambda,x,t)q(t)y(\lambda,t)dt=i\int\limits_0^xK_1(\lambda,x,t)q(t)y_0(\lambda,t)dt}\\
{\displaystyle+i\int\limits_0^xK_1(\lambda,x,t)q(t)z(\lambda,t)dt}
\end{array}\label{eq2.5}
\end{equation}
where
\begin{equation}
K_1(\lambda,x,t)\stackrel{\rm def}{=}\frac{s_2(i\lambda(x-t))}{(i\lambda)^2}.\label{eq2.6}
\end{equation}
In $L^2(0,l)$, define the family of Volterra operators
\begin{equation}
(K_\lambda f)(x)\stackrel{\rm def}{=}\int\limits_0^xK_1(\lambda,x,t)q(t)f(t)dt\quad(f\in L^2(0,l))\label{eq2.7}
\end{equation}
depending on $\lambda\in\mathbb{C}$, then equation \eqref{eq2.5}, in terms of $K_\lambda$, becomes
$$(I-iK_\lambda)z(\lambda,x)=iK_\lambda y_0(\lambda,x),$$
and thus
\begin{equation}
z(\lambda,x)=\sum\limits_{n=1}^\infty i^nK_\lambda^ny_0(\lambda,x).\label{eq2.8}
\end{equation}
The operators $K_\lambda^n$ are Volterra also,
$$(K_\lambda^nf)(x)\stackrel{\rm def}{=}\int\limits_0^xK_n(\lambda,x,t)q(t)f(t)dt,$$
and, for the kernels $K_n(\lambda,x,t)$, the following recurrent relations hold:
\begin{equation}
K_{n+1}(\lambda,x,t)=\int\limits_t^xK_1(\lambda,x,s)q(s)K_n(\lambda,s,t)ds\quad(n\in\mathbb{N}).\label{eq2.9}
\end{equation}
Using $\lambda=\alpha+i\beta\in\mathbb{C}$ ($\alpha$, $\beta\in\mathbb{R}$), we obtain
$$i\lambda=-\beta+i\alpha,\quad i\lambda\zeta_2=\frac12(\beta-\alpha\sqrt3)-\frac i2(\alpha+\beta\sqrt3);\quad i\lambda\zeta_3=\frac12(\beta+\alpha\sqrt3)-\frac i2(\alpha-\beta\sqrt3);$$
and thus
$$s_p(i\lambda)=\frac13\left\{e^{-\beta+i\alpha}+\frac1{\zeta_2^p}e^{\frac12(\beta-\alpha\sqrt3)-\frac i2(\alpha+\beta\sqrt3)}+\frac1{\zeta_3^p}e^{\frac12(\beta+\alpha\sqrt3)-\frac i2(\alpha-\beta\sqrt3)}\right\}\quad(p=0,1,2),$$
therefore
$$|s_p(i\lambda)|\leq\frac13\left(e^{-\beta}+2e^\frac\beta2\ch\frac{\alpha\sqrt3}2\right).$$
Hence it follows that
\begin{equation}
|s_p(i\lambda)|\leq d(\lambda)\quad(d(\lambda)\stackrel{\rm def}{=}e^{|\beta|}\ch\frac{\alpha\sqrt3}2),\label{eq2.10}
\end{equation}
thus, for $K_1(\lambda,x,t)$ \eqref{eq2.6}, the following estimate holds:
\begin{equation}
|K_1(\lambda,x,t)|\leq\frac1{|\lambda|^2}d(\lambda(x-t))\quad(\lambda\not=0).\label{eq2.11}
\end{equation}

\begin{lemma}\label{l2.1}
The kernels $K_n(\lambda,x,t)$ \eqref{eq2.9} have the properties
\begin{equation}
K_n(\lambda\zeta_2,x,t)=K_n(\lambda,x,t);\quad\overline{K_n(\lambda,x,t)}=K_n(\overline{\lambda},t,x);\label{eq2.12}
\end{equation}
and satisfy the inequalities
\begin{equation}
|K_n(\lambda,x,t)|\leq\left\{
\begin{array}{lll}
{\displaystyle\frac{d(\lambda(x-t))}{|\lambda|^{2n}}\frac{\sigma^{n-1}(x)}{(n-1)!}\quad(\lambda\not=0);}\\
{\displaystyle\left(\frac{(x-t)^2}2\right)^n\frac{\sigma^{n-1}(x)}{n^{2n}(n-1)!}\quad(\lambda=0)}
\end{array}\quad(n\in\mathbb{N})\right.\label{eq2.13}
\end{equation}
where $d(\lambda)$ is given by \eqref{eq2.10} ($\lambda=\alpha+i\beta$; $\alpha\in\mathbb{R}$; $\beta\in\mathbb{R}$), and
\begin{equation}
\sigma(x)\stackrel{\rm def}{=}\int\limits_0^x|q(t)|dt.\label{eq2.14}
\end{equation}
\end{lemma}

P r o o f. Relations (ii), (iii) \eqref{eq1.4} imply that the kernel $K_1(\lambda,x,t)$ \eqref{eq2.1} satisfies relation \eqref{eq2.12}, whence, due to \eqref{eq2.9}, follows that $K_n(\lambda,x,t)$ also has properties \eqref{eq2.12}, because $q(x)$ is real.

Consider the case of $\lambda\not=0$. Inequality \eqref{eq2.13} ($\lambda\not=0$), for $n=1$, coincides with \eqref{eq2.11}. Using induction on $n$ and recurrent relations \eqref{eq2.9}, and also the estimates \eqref{eq2.11} and \eqref{eq2.13} (for $n$), we obtain
$$|K_{n+1}(\lambda,x,t)|\leq\int\limits_t^x\frac1{|\lambda|^2}d(\lambda(x-s))|q(s)|\frac{d(\lambda(s-t))}{|\lambda|^{2n}}\frac{\sigma^{n-1}(s)}{(n-1)!}ds.$$
And since ($\lambda=\alpha+i\beta$; $\alpha$, $\beta\in\mathbb{R}$)
$$d(\lambda(x-s))d(\lambda(s-t))=e^{|\beta|(x-t)}\ch\frac{\alpha\sqrt3}2(x-s)\ch\frac{\alpha\sqrt3}2(s-t)=\frac12e^{|\beta|(x-t)}$$
$$\times\left[\ch\frac{\alpha\sqrt3}2(x-t)+\ch\frac{\alpha\sqrt3}2(x+t-2s)\right]\leq d(\lambda(x-t)),$$
then
$$|K_{n+1}(\lambda,x,t)|\leq\frac{d(\lambda(x-t))}{|\lambda|^{2(n+1)}}\frac{\sigma^n(x)}{n!},$$
which proves \eqref{eq2.13} for $n+1$ ($\lambda\not=0$).

Prove inequality \eqref{eq2.13} when $\lambda=0$. Equation \eqref{eq2.6} yields that $|K_1(0,x,t)|=$ ${\displaystyle\frac{(x-t)^2}2}$, which coincides with \eqref{eq2.13} for $\lambda=0$. Again using induction on $n$ and \eqref{eq2.9}, we obtain
$$|K_{n+1}(0,x,t)|\leq\int\limits_t^x\frac{(x-s)^2}2|q(s)|\frac{(s-t)^{2n}}{2^n}\frac{\sigma^{n-1}(s)}{n^{2n}(n-1)!}ds.$$
Function $f(s)=(x-s)^2(s-t)^{2n}$ is positive for $s\in(t,x)$ and $f(t)=f(s)=0$. On the interval $(t,x)$, it reaches its maximum value  at the point $s_0=(t+nx)/(1+n)$ and ${\displaystyle f(s_0)=\frac{(x-t)^{2(n+1)}}{(1+n)^{2(n+1)}}\cdot n^{2n}}$. Hence it follows that
$$|K_{n+1}(0,x,t)|\leq\frac{(x-t)^{2(n+1)}}{2^{n+1}}\cdot\frac1{(n+1)^{2(n+1)}}\cdot\frac{\sigma^n(x)}{n!}.\blacksquare$$

Rewrite equality \eqref{eq2.8} as
\begin{equation}
z(\lambda,x)=\int\limits_0^xN(\lambda,x,t)q(t)y_0(\lambda,t)dt\label{eq2.15}
\end{equation}
where
\begin{equation}
N(\lambda,x,t)\stackrel{\rm def}{=}\sum\limits_{n=1}^\infty i^nK_n(\lambda,x,t).\label{eq2.16}
\end{equation}
Convergence of series \eqref{eq2.16} follows from \eqref{eq2.13} and
\begin{equation}
|N(\lambda,x,t)|\leq\left\{
\begin{array}{lll}
{\displaystyle\frac{d(\lambda(x-t))}{|\lambda|^2}\cdot\exp\left\{\frac{\sigma(x)}{|\lambda|^2}\right\}\quad(\lambda\not=0);}\\
{\displaystyle\frac{(x-t)^2}2\cdot\exp\left\{\frac{(x-t)^2\sigma(x)}2\right\}\quad(\lambda=0).}
\end{array}\right.\label{eq2.17}
\end{equation}

\begin{lemma}\label{l2.2}
Solution $y(\lambda,x)$ to the Cauchy problem \eqref{eq2.2}, \eqref{eq2.3} is expressed via the solution $y_0(\lambda,x)$ to the Cauchy problem \eqref{eq1.5} ($y_0=0$, $f=0$) by the formula
\begin{equation}
y(\lambda,x)=(I+T_\lambda)y_0(\lambda,x)\label{eq2.18}
\end{equation}
where $T_\lambda$ is a family of Volterra operators in $L^2(0,l)$,
\begin{equation}
\begin{array}{ccc}
{\displaystyle(T_\lambda f)(x)\stackrel{\rm def}{=}\int\limits_0^xT(\lambda,x,t)f(t)dt;\quad T(\lambda,x,t)\stackrel{\rm def}{=}N(\lambda,x,t)q(t);}\\ T(\lambda\zeta_2,x,t)=T(\lambda,x,t),
\end{array}\label{eq2.19}
\end{equation}
and $N(\lambda,x,t)$ is given by \eqref{eq2.16}, $f\in L^2(0,l)$.

Operators $K_\lambda$ \eqref{eq2.7} and $T_\lambda$ \eqref{eq2.19} satisfy the identities
\begin{equation}
T_\lambda(I-iK_\lambda)^{-1}iK_\lambda=(I-iK_\lambda)^{-1}-I;\quad I+T_\lambda=(I-iK_\lambda)^{-1}.\label{eq2.20}
\end{equation}
The kernel $T(\lambda,x,t)$ is the solution to the integral equation
\begin{equation}
T(\lambda,x,t)-i\int\limits_t^xK_1(\lambda,x,s)q(s)T_1(\lambda,s,t)ds=iK_1(\lambda,x,t)q(t)\label{eq2.21}
\end{equation}
where $K_1(\lambda,x,t)$ is given by \eqref{eq2.6} and
\begin{equation}
\lim\limits_{t\rightarrow x}\frac{T(\lambda,x,t)}{(x-t)^2}=\frac i2q(x).\label{eq2.22}
\end{equation}
\end{lemma}

P r o o f. We need to prove the equalities \eqref{eq2.21} and \eqref{eq2.22}. The equation $(I-iK_\lambda)(I+T_\lambda)=I$ \eqref{eq2.20} implies that $T_\lambda-iK_\lambda-iK_\lambda T_\lambda=0$, or
$$\int\limits_0^xT(\lambda,x,t)f(t)dt-i\int\limits_0^xK_1(\lambda,x,t)q(t)f(t)dt$$
$$-i\int\limits_0^xdtK_1(\lambda,x,t)q(t)\int\limits_0^tT(\lambda,t,s)f(s)ds=0,$$
i. e.,
$$\int\limits_0^xdtf(t)\left\{T(\lambda,x,t)-iK_1(\lambda,x,t)q(t)-i\int\limits_t^sK_1(\lambda,x,s)q(s)T(\lambda,s,t)ds\right\}=0,$$
whence, due to arbitrariness of $f\in L^2(0,l)$, follows \eqref{eq2.21}. Equality \eqref{eq2.22} is a corollary of relation \eqref{eq2.21}, in accordance with (xi) \eqref{eq1.4}. $\blacksquare$

Operators $I+T_\lambda$ \eqref{eq2.18} are analogues of {\bf transformation operators} \cite{1,2,3} of the pair of operators $\{L_0(\theta),L_q(\theta)\}$.
\vspace{5mm}

{\bf 2.2} Define the functions $s_p(\lambda,x)$ obtained from $s_p(i\lambda x)/(i\lambda)^p$ using the transformation operators $I+T_\lambda$ \eqref{eq2.18}
\begin{equation}
s_p(\lambda,x)\stackrel{\rm def}{=}\frac{s_p(i\lambda x)}{(i\lambda)^p}+\int\limits_0^xT(\lambda,x,t)\frac{s_p(i\lambda t)}{(i\lambda)^p}dt\quad(0\leq p\leq2)\label{eq2.23}
\end{equation}
which are solutions to equation \eqref{eq2.2} and satisfy the initial data
$$s_0(\lambda,0)=1;\quad s'_0(\lambda,0)=0;\quad s''_0(\lambda,0)=0;$$
\begin{equation}
s_1(\lambda,0)=0;\quad s'_1(\lambda,0)=1;\quad s''_1(\lambda,0)=0;\label{eq2.24}
\end{equation}
$$s_2(\lambda,0)=0;\quad s'_2(\lambda,0)=0;\quad s''_2(\lambda,0)=1.$$
Moreover,
\begin{equation}
s_p(\lambda\zeta_2,x)=s_p(\lambda,x)\quad(0\leq p\leq2),\label{eq2.25}
\end{equation}
due to (ii) \eqref{eq1.4} and \eqref{eq2.19}. Rewrite equality \eqref{eq2.23} as
$$3(i\lambda)^ps_p(\lambda,x)e^{-i\lambda\zeta_kx}=3s_p(i\lambda x)e^{-i\lambda\zeta_kx}+\int\limits_0^xT(\lambda,x,t)e^{-i\lambda\zeta_k(x-t)}3s_p(i\lambda t)e^{-i\lambda\zeta_kt}dt,$$
then, taking \eqref{eq1.16} into account, we find asymptotic of $\{s_p(\lambda,x)\}$ in the sectors $\{\Omega_k\}$ \eqref{eq1.14}.

\begin{lemma}\label{l2.3}
For entire functions of exponential type $\{s_p(\lambda,x)\}_0^2$ \eqref{eq2.23}, the following relations are true:
\begin{equation}
\lim\limits_{\lambda\rightarrow\infty}3\zeta_k^p(i\lambda)^ps_p(\lambda,x)e^{-i\lambda\zeta_kx}=1\quad(\lambda\in\Omega_k)\label{eq2.26}
\end{equation}
($0\leq p\leq2$, $1\leq k\leq3$).
\end{lemma}

Function
\begin{equation}
Y_0(\lambda,x)\stackrel{\rm def}{=}y_1s_1(\lambda,x)+y_2s_2(\lambda,x)\quad(y_1,y_2\in\mathbb{C})\label{eq2.27}
\end{equation}
is the solution to equation \eqref{eq2.2} and $Y_0(\lambda,0)=0.$ The second and the third boundary conditions in \eqref{eq1.18} for $Y_0(\lambda,x)$ lead to the system of linear equations
\begin{equation}
\left\{
\begin{array}{lll}
y_1(s'_1(\lambda,l)-\theta)+y_2s'_2(\lambda,l)=0;\\
y_1s_1(\lambda,l)+y_2s_2(\lambda,l)=0;
\end{array}\right.\label{eq2.28}
\end{equation}
coinciding with \eqref{eq1.20} when $q\equiv0$. This system has a non-zero solution $y_1$, $y_2$ only if its determinant $\Delta_\theta(q,\lambda)$ vanishes, $\Delta_\theta(q,\lambda)=0$, where
\begin{equation}
\Delta_\theta(q,\lambda)\stackrel{\rm def}{=}s_2(\lambda,l)s'_1(\lambda,l)-s_1(\lambda,l)s'_2(\lambda,l)-\theta s_2(\lambda,l).\label{eq2.29}
\end{equation}
By $W_{k,p}(\lambda,x)$, we denote the Wronskian,
\begin{equation}
W_{k,p}(\lambda,x)\stackrel{\rm def}{=}s_k(\lambda,x)s'_p(\lambda,x)-s_p(\lambda,x)s'_k(\lambda,x)\quad(0\leq k,p\leq2),\label{eq2.30}
\end{equation}
and define the {\bf operation} ``*'',
\begin{equation}
f^*(\lambda)=\overline{f(\overline{\lambda})}.\label{eq2.31}
\end{equation}

\begin{lemma}\label{l2.4}
For the Wronskians $\{W_{k,s}(\lambda,x)\}$ \eqref{eq2.30}, the following representations are true:
\begin{equation}
W_{0,1}(\lambda,x)=s_0^*(\lambda,x);\quad W_{1,2}(\lambda,x)=s_2^*(\lambda,x);\quad W_{0,2}(\lambda,x)=s_1^*(\lambda,x)\label{eq2.32}
\end{equation}
where $\{s_p(\lambda,x)\}$ are given by \eqref{eq2.23}.
\end{lemma}

P r o o f. Derivative of the function
$$W_{0,1}(\lambda,x)=s_0(\lambda,x)s'_1(\lambda,x)-s_1(\lambda,x)s'_0(\lambda,x)\quad(W_{0,1}(\lambda,0)=1)$$
equals
$$W'_{0,1}(\lambda,x)=s_0(\lambda,x)s''_1(\lambda,x)-s_1(\lambda,x)s''_0(\lambda,x)\quad(W'_{0,1}(\lambda,0)=0,$$
and thus
$$W''_{0,1}(\lambda,x)=s'_0(\lambda,x)s''_1(\lambda,x)-s'_1(\lambda,0)s''_0(\lambda,x)\quad(W''_{0,1}(\lambda,0)=0).$$
Hence, due to \eqref{eq2.2}, it follows that $W_{0,1}(\lambda,x)$ is a solution to the equation
$$iy'''(x)-q(x)y(x)=-\lambda^3y(x)$$
and satisfies the initial data $y(0)=1$, $y'(0)=0$, $y''(0)=0$. This equation is derived from \eqref{eq2.2} upon the complex conjugation, and taking into account \eqref{eq2.24}, we obtain that $W_{0,1}(\lambda,x)=s_0^*(\lambda,x)$. Other relations in \eqref{eq2.32} are proved analogously. $\blacksquare$

For the functions $s_p^*(\lambda,x)$, analogously to \eqref{eq2.23}, the following relations hold:
\begin{equation}
s_p^*(\lambda,x)=\frac{s_p^*(-i\lambda x)}{(-i\lambda)^p}+\int\limits_0^xT^*(\lambda,x,t)\frac{s_p^*(-i\lambda t)}{(-i\lambda)^p}dt\quad(0\leq p\leq2).\label{eq2.33}
\end{equation}

\begin{lemma}\label{l2.5}
Characteristic function $\Delta_\theta(q,\lambda)$ \eqref{eq2.29} equals
\begin{equation}
\Delta_\theta(q,\lambda)=-\{\theta s_2(\lambda,l)+s_2^*(\lambda,l)\}\label{eq2.34}
\end{equation}
where $s_2(\lambda,x)$ and $s_2^*(\lambda,x)$ are given by \eqref{eq2.23} and \eqref{eq2.33}, besides
\begin{equation}
\Delta_\theta(q,\lambda\zeta_2)=\Delta_\theta(q,\lambda);\quad\overline{\Delta_\theta(q,\lambda)}=\overline{\theta}\Delta_\theta(q,\overline{\lambda}).\label{eq2.35}
\end{equation}
\end{lemma}

The function $\Delta_\theta(q,\lambda)$ \eqref{eq2.34} coincides with $\Delta_\theta(0,\lambda)$ \eqref{eq1.21} for $q\equiv0$. Every root $\lambda_n(q,\theta)$ of the equation $\Delta_\theta(q,\lambda)=0$ is included with the series $\{\zeta_2^p\lambda_n(q,\theta)\}$ ($0\leq p\leq2$), due to \eqref{eq2.35}. Enumerate real zeros $\lambda_n(q,\theta)\in\mathbb{R}$ of the function $\Delta_\theta(q,\lambda)$ in ascending order $...\lambda_{-1}(q,\theta)<\lambda_0(q,\theta)<\lambda_1(q,\theta)<...$, here $\lambda_0(q,\theta)$ is the smallest non-negative zero of $\Delta_\theta(q,\lambda)$.

Equations \eqref{eq2.23} and \eqref{eq2.10}, \eqref{eq2.17} imply that $s_2(\lambda,x)$ is an entire function of exponential type (for all $x\in[0,l]$),
\begin{equation}
\begin{array}{cccc}
{\displaystyle|s_2(\lambda,x)|\leq\frac{d(\lambda x)}{|\lambda|^2}+\int\limits_0^x|q(t)|\frac{d(\lambda(x-t))}{|\lambda|^2}\exp\left\{\frac{\sigma(x)}{|\lambda|^2}\right\}\frac{dt}{|\lambda|^2}}\\
{\displaystyle\leq\frac{d(\lambda x)}{|\lambda|^2}\left\{1+\frac{\sigma(x)}{|\lambda|^2}\cdot\exp\left(\frac{\sigma(x)}{|\lambda|^2}\right)\right\}.}
\end{array}\label{eq2.36}
\end{equation}
Therefore $\Delta_\theta(q,\lambda)$ \eqref{eq2.34} also is an entire function of exponential type and for it the multiplicative expansion \eqref{eq1.26} (Theorem \ref{t1.1}) is true. Product of primary factors of the series $\{\zeta_2^l\lambda_n(q,\theta)\}$ ($l=0$, 1, 2) of a real root $\lambda_n(q,\theta)$ equals $1-\lambda^3/\lambda_n^3(q,\theta)$, and thus
$$\Delta_\theta(q,\lambda)=\alpha e^{\beta\lambda}\lambda^{3m}\prod\limits_n\left(1-\frac{\lambda^3}{\lambda_n^3(q,\theta)}\right)$$
where $m\in\mathbb{Z}_+$; $\alpha$, $\beta\in\mathbb{C}$. Equation \eqref{eq2.34} yields that $\Delta_\theta(q,0)=-(\theta\alpha+\overline{\alpha}))$, here
$$\alpha\stackrel{\rm def}{=}s_2(0,l)=\frac{l^2}2+\int\limits_0^lT(0,l,t)\frac{t^2}2dt(\not=0),$$
therefore $\theta_0\stackrel{\rm def}{=}\overline{a}/a\in\mathbb{T}$ and $\Delta_\theta(q,0)\not=0$ for $\theta+\theta_0\not=0$ (in the case of $q(x)\equiv0$, $\theta_0=1$, see Lemma \ref{l1.5}). $\beta=0$ since ${\displaystyle\left.\frac d{d\lambda}(\Delta_\theta(q,\lambda))\right|_{\lambda=0}=0}$ ($\theta+\theta_0\not=0$).

\begin{lemma}\label{l2.6}
Characteristic function $\Delta_\theta(q,\lambda)$ \eqref{eq2.34} has a multiplicative representation:
\begin{equation}
\Delta_\theta(q,\lambda)=-a(\theta+\theta_0)\prod\limits_n\left(1-\frac{\lambda^3}{\lambda_n^3(q,\theta)}\right)\quad(\theta+\theta_0\not=0)\label{eq2.37}
\end{equation}
where $a=s_2(0,l)\not=0$ depends only on $q(x)$; $\theta_0=\overline{a}/a\in\mathbb{T}$; $\lambda_n(q,\theta)$ are the real zeros of $\Delta_\theta(q,\lambda)$ enumerated in ascending order.
\end{lemma}

If $\theta+\theta_0=0$, then in representation \eqref{eq2.37} of the function $\Delta_\theta(q,\lambda)$ the factor $\lambda^3$ appears, and thus the operator $L_q(\theta)$ has the eigenvalue $\lambda^3=0$. We confine ourselves to the case of $\theta+\theta_0\not=0$.
\vspace{5mm}

{\bf 2.3} Study the asymptotic behavior of zeros $\lambda_n(q,\theta)$ when $|n|\rightarrow\infty$. Substitute the expressions \eqref{eq2.23}, for $s_2(\lambda,x)$, and \eqref{eq2.33}, for $s_2^*(\lambda,x)$, into formula \eqref{eq2.34}, then, taking into account \eqref{eq1.21}, we obtain
\begin{equation}
\Delta_\theta(q,\lambda)=\Delta_\theta(0,\lambda)+Q_\theta(\lambda),\label{eq2.38}
\end{equation}
here
\begin{equation}
Q_\theta(\lambda)\stackrel{\rm def}{=}-\int\limits_0^l\left\{\theta T(\lambda,l,t)\frac{s_2(i\lambda t)}{(i\lambda)^2}+T^*(\lambda,l,t)\frac{s_2(-i\lambda t)}{(i\lambda)^2}\right\}dt.\label{eq2.39}
\end{equation}

\begin{lemma}\label{l2.7}
For the function $Q_\theta(\lambda)$ \eqref{eq2.39}, the following representation is true:
\begin{equation}
Q_0(\lambda)=-i\int\limits_0^lq(t)\left\{\theta s_2(\lambda,t)\frac{s_2(i\lambda(l-t))}{(i\lambda)^2}-s_2^*(\lambda,t)\frac{s_2(-i\lambda(l-t))}{(-i\lambda)^2}\right\}dt\label{eq2.40}
\end{equation}
where $s_2(\lambda,x)$ and $s_2^*(\lambda,x)$ are given by \eqref{eq2.23} and \eqref{eq2.33}.
\end{lemma}

P r o o f. Relations \eqref{eq2.16}, \eqref{eq2.19} imply
$$A\stackrel{\rm def}{=}\int\limits_0^lT(\lambda,l,t)s_2(i\lambda t)dt=\sum\limits_1^\infty i^n\int\limits_0^lK_n(\lambda,l,t)s_2(i\lambda t)dt,$$
and, due to \eqref{eq2.6}, \eqref{eq2.9},
$$A=\frac i{(i\lambda)^2}\int\limits_0^lq(t)s_2(i\lambda(l-t))s_2(i\lambda t)dt+\frac{i^2}{(i\lambda)^4}\int\limits_0^l\int\limits_t^ls_2(i\lambda(l-s))q(s)s_2(i\lambda(s-t))ds$$
$$\times q(t)s_2(i\lambda t)dt+...=\frac i{(i\lambda)^2}\int\limits_0^lq(t)s_2(i\lambda(l-t))s_2(i\lambda t)dt$$
$$+\frac{(i)^2}{(i\lambda)^4}\int\limits_0^ls_2(i\lambda(l-s))q(s)\int\limits_0^ss_2(i\lambda(s-t))q(t)s_2(i\lambda t)dt+...=\frac i{(i\lambda)^2}\int\limits_0^lq(t)s_2(i\lambda(l-t))$$
$$\times(1-iK_\lambda)^{-1}s_2(i\lambda t)dt.$$
Using \eqref{eq2.20}, we obtain
$$A=\frac i{(i\lambda)^2}\int\limits_0^lq(t)s_2(i\lambda(l-t))(I+T_\lambda)s_2(i\lambda t)=i\int\limits_0^lq(t)s_2(i\lambda(l-t))s_2(\lambda,t).$$
Analogously, it is proved that
$$B\stackrel{\rm def}{=}\int\limits_0^lT^*(\lambda,l,t)s_2(-i\lambda t)dt=-i\int\limits_0^lq(t)s_2(-i\lambda(l-t))s_2^*(\lambda,t)dt.\blacksquare$$

Taking into account \eqref{eq2.36}, for $Q_\theta(\lambda)$ \eqref{eq2.40} we find
\begin{equation}
\begin{array}{ccc}
{\displaystyle|Q_\theta(\lambda)|\leq\frac2{|\lambda|^2}\int\limits_0^l|q(t)|\frac{d(\lambda t)}{|\lambda|^2}\left(1+\frac{\sigma(t)}{|\lambda|^2}\exp\left\{\frac{\sigma(t)}{|\lambda|^2}\right\}\right)d(\lambda(l-t))dt}\\
{\displaystyle\leq\frac2{|\lambda|^4}d(\lambda l)\left(\sigma(l)+\frac{\sigma^2(l)}{2|\lambda|^2}\exp\left\{\frac{\sigma(l)}{|\lambda|^2}\right\}\right).}
\end{array}\label{eq2.41}
\end{equation}
Relation \eqref{eq2.10} implies that for the characteristic function $\Delta_\theta(0,\lambda)$ \eqref{eq1.21} the following estimate is true:
\begin{equation}
|\Delta_\theta(0,\lambda)|\leq\frac{2d(\lambda l)}{|\lambda|^2}\label{eq2.42}
\end{equation}
where $d(\lambda)$ is given by \eqref{eq2.10}. Thus, due to Rouche's theorem, the functions $\Delta_\theta(0,\lambda)$ and $\Delta_\theta(0,\lambda)+Q_\theta(\lambda)(=\Delta_\theta(q,\lambda))$ inside the circle $|\lambda|<R$ ($R\gg1$) have the same number of zeros. Moreover, using \eqref{eq2.41}, \eqref{eq2.42}, again according to Rouche's theorem, in the neighborhood $|\lambda_n(0,\theta)-\lambda|<r_n$ ($0<r_n\ll1$, $|\lambda_n(0,\theta)|>R$ ($\gg1$)) there lies exactly one zero $\lambda_n(q,\theta)$ of the function $\Delta_\theta(q,\lambda)$, and thus
\begin{equation}
\lambda_n(q,\theta)=\lambda_n(0,\theta)+o\left(\frac1{\lambda_n^2(0,\theta)}\right).\label{eq2.43}
\end{equation}

\begin{theorem}\label{t2.1}
Spectrum $\sigma(L_q)$ of the operator $L_q(\theta)$ \eqref{eq2.1} equals
\begin{equation}
\sigma(L_q)\stackrel{\rm def}{=}\{\lambda_n^3(q,\theta):n\in\mathbb{Z}\}\label{eq2.44}
\end{equation}
where $\{\lambda_n(q,\theta)\}$ are real zeros of characteristic function $\Delta_\theta(q,\lambda)$ \eqref{eq2.34} enumerated in ascending order and having, for $|n|\rightarrow\infty$, asymptotic \eqref{eq2.43}.

Eigenfunctions of $L_q(\theta)$ corresponding to $\lambda=\lambda_n(q,\theta)$ are given by
\begin{equation}
\psi_n(q,\lambda,x)=\frac1{a_n(\lambda)}\{s_2(\lambda,x)s_1(\lambda,l)-s_1(\lambda,x)s_2(\lambda,l)\},\label{eq2.45}
\end{equation}
besides, $\psi_n(q,\lambda\zeta_2,x)=\psi_n(q,\lambda,x)$ and the numbers $a_n(\lambda)$ ($>0$) are defined from the condition $\|\psi_n(q,\lambda,x)\|_{L^2}=1$.
\end{theorem}

To obtain an analogue of representation \eqref{eq1.31} for $\psi_n(q,\lambda,x)$ \eqref{eq2.45}, substitute expressions \eqref{eq2.23} into the formula
$$B\stackrel{\rm def}{=}s_2(\lambda,x)s_1(\lambda,t)-s_1(\lambda,x)s_2(\lambda,t)=\frac1{(i\lambda)^3}\{s_2(i\lambda x)s_1(i\lambda t)-s_1(i\lambda x)s_2(i\lambda x)$$
$$+\int\limits_0^xd\tau T(\lambda,x,\tau)[s_2(i\lambda\tau)s_1(i\lambda t)-s_1(i\lambda\tau)s_2(i\lambda t)]$$
$$+\int\limits_0^td\eta T(\lambda,t,\eta)[s_2(i\lambda x)s_1(i\lambda\eta)-s_1(i\lambda x)s_2(i\lambda\eta)]$$
$$\left.+\int\limits_0^xd\tau T(\lambda,x,\tau)\int\limits_0^td\eta T(\lambda,t,\eta)[s_2(i\lambda\tau)s_1(i\lambda\eta)-s_1(i\lambda\tau)s_2(i\lambda\eta)]\right\}.$$
Using \eqref{eq1.30}, we obtain that
$$B=(I+T_\lambda)_x(I+T_\lambda)_t\frac1{\sqrt3\lambda^3}\{s_0(i\lambda(x+\zeta_2t))-s_0(i\lambda(x+\zeta_3t))\}$$
where operators $(I+T_\lambda)_x$ and $(I+T_\lambda)_t$ are given by \eqref{eq2.18}, \eqref{eq2.19} and act with respect to the variables $x$ and $t$ correspondingly. Operators $(I+T_\lambda)_x$ and $(I+T_\lambda)_t$ commute.

\begin{remark}\label{r2.1}
The eigenfunctions $\psi_n(q,\lambda,x)$ \eqref{eq2.45} have the representation
$$\psi_n(q,\lambda,x)=\frac1{\sqrt3\lambda^2a_n(\lambda)}\left.\{(I+T_\lambda)_x(I+T_\lambda)_t[s_0(i\lambda(x+\zeta_2t))-s_0(i\lambda(x+\zeta_3t))]\}\right|_{tl},$$
besides, $\lambda=\lambda_n(q,\theta)$.
\end{remark}
\vspace{5mm}

{\bf 2.4} Calculate the resolvent $R_{L_q}(\lambda^3)=(L_q(\theta)-\lambda^3I)^{-1}$. Consider the Cauchy problem
\begin{equation}
iD^3y(x)+q(x)y(x)=\lambda^3y(x)+f(x)\,(x\in\mathbb{R}_+);\quad y(0)=y_0,\, y'(0)=y_1,\, y''(0)=y_2.\label{eq2.46}
\end{equation}
The function
\begin{equation}
Y_0(\lambda,x)=y_0s_0(\lambda,x)+y_1s_1(\lambda,x)+y_2s_2(\lambda,x),\label{eq2.47}
\end{equation}
where $\{s_p(\lambda,x)\}_0^2$ are given by \eqref{eq2.23}, is the solution to Cauchy problem \eqref{eq2.46} when $f=0$. Using the method of variation of constants, we find the general solution to problem \eqref{eq2.46}:
\begin{equation}
Y(\lambda,x)=Y_0(\lambda,x)-i\int\limits_0^xG(\lambda,x,t)f(t)dt\label{eq2.48}
\end{equation}
where $Y_0(\lambda,x)$ is given by \eqref{eq2.47} and the kernel $G(\lambda,x,t)$ equals
\begin{equation}
G(\lambda,x,t)\stackrel{\rm def}{=}\det\left[
\begin{array}{ccc}
s_0(\lambda,t)&s_1(\lambda,t)&s_2(\lambda,t)\\
s'_0(\lambda,t)&s'_1(\lambda,t)&s'_2(\lambda,t)\\
s_0(\lambda,x)&s_1(\lambda,x)&s_2(\lambda,x)
\end{array}\right].\label{eq2.49}
\end{equation}

\begin{lemma}\label{l2.8}
The function $G(\lambda,x,t)$ \eqref{eq2.49} is the solution to the equation
\begin{equation}
i\frac{d^3}{dx^3}G(\lambda,x,t)+q(x)G(\lambda,x,t)=\lambda^3G(\lambda,x,t)\label{eq2.50}
\end{equation}
and
\begin{equation}
\left.G(\lambda,x,t)\right|_{x=t}=0;\quad\left.\frac d{dx}G(\lambda,x,t)\right|_{x=t}=0;\quad\left.\frac{d^2}{dx^2}G(\lambda,x,t)\right|_{x=t}=1.\label{eq2.51}
\end{equation}
The following representation is true:
\begin{equation}
G(\lambda,x,t)=s_0(\lambda,x)s_2^*(\lambda,t)-s_1(\lambda,x)s_1^*(\lambda,t)+s_2(\lambda,x)s_0^*(\lambda,t),\label{eq2.52}
\end{equation}
here $\{s_p(\lambda,x)\}$ and $\{s_p^*(\lambda,x)\}$ are given by \eqref{eq2.23} and \eqref{eq2.33}, besides, $\overline{G(\lambda,x,t)}=G(\overline{\lambda},t,x)$.
\end{lemma}

P r o o f. The function
$$F(\lambda,x)\stackrel{\rm def}{=}\left.\frac{d^2}{dx^2}G(\lambda,x,t)\right|_{t=x}=\det\left[
\begin{array}{ccc}
s_0(\lambda,x)&s_1(\lambda,x)&s_2(\lambda,x)\\
s'_0(\lambda,x)&s'_1(\lambda,x)&s'_2(\lambda,x)\\
s''_0(\lambda,x)&s''_1(\lambda,x)&s''_2(\lambda,x)
\end{array}\right]$$
does not depend on $x$ because $F'(\lambda,x)=0$ and, taking into account \eqref{eq2.24}, we obtain $F(\lambda,x)=1$, which gives the last equality in \eqref{eq2.51}. Equation \eqref{eq2.49} implies that
$$G(\lambda,x,t)=s_0(\lambda,x)W_{1,2}(\lambda,t)-s_1(\lambda,x)W_{0,2}(\lambda,t)+s_2(\lambda,x)W_{0,1}(\lambda,t),$$
which, in view of \eqref{eq2.32}, gives \eqref{eq2.52}. The equality $\overline{G(\lambda,x,t)}=G(\overline{\lambda},t,x)$ follows from representation \eqref{eq2.52}. $\blacksquare$

Calculate the resolvent $R_{L_q}(\lambda^3)$, and let $Y=R_{L_q}(\lambda^3)f$, then $L_q(\theta)Y-\lambda^3Y=f$, and thus $Y=Y(\lambda,x)$ is the solution to equation \eqref{eq2.46}, moreover,
\begin{equation}
Y(\lambda,x)=y_1s_1(\lambda,x)+y_2s_2(\lambda,x)-i\int\limits_0^xG(\lambda,x,t)f(t)dt\label{eq2.53}
\end{equation}
and $Y(\lambda,0)=0$. The second and the third boundary conditions in \eqref{eq1.18} imply the following system of equations:
$$\left\{
\begin{array}{lll}
{\displaystyle y_1(s'_1(\lambda,l)-\theta)+y_2s'_2(\lambda,l)=i\int\limits_0^lG'(\lambda,l,t)f(t)dt;}\\
{\displaystyle y_1s_1(\lambda,l)+y_2s_2(\lambda,l)=i\int\limits_0^lG(\lambda,l,t)f(t)dt}
\end{array}\right.$$
(${\displaystyle G'(\lambda,l,t)=\left.\frac d{dx}G(\lambda,x,t)\right|_{x=l}}$), determinant of which equals $\Delta_\theta(q,\lambda)$ \eqref{eq2.23}. Hence, for $\Delta_\theta(q,\lambda)\not=0,$ we find
$$y_1=\frac i{\Delta_\theta(q,\lambda)}\int\limits_0^l[G'(\lambda,l,t)s_2(\lambda,l)-G(\lambda,l,t)s'_2(\lambda,l)]f(t)dt;$$
$$y_2=\frac i{\Delta_\theta(q,\lambda)}\int\limits_0^l[G(\lambda,l,t)(s'_1(\lambda,l)-\theta)-G'(\lambda,l,t)s_1(\lambda,l)]f(t)dt,$$
and upon substituting these expressions into \eqref{eq2.53}, we obtain that
$$Y(\lambda,x)=\frac i{\Delta_\theta(q,\lambda)}\left\{\int\limits_0^l[s_1(\lambda,x)(G'(\lambda,l,t)s_2(\lambda,l)-G(\lambda,l)s'_2(\lambda,lt))\right.$$
\begin{equation}
+s_2(\lambda,x)(G(\lambda,l,t)s'_1(\lambda,l)-G'(\lambda,l,t)s_1(\lambda,l)-\theta G(\lambda,l,t))]f(t)dt\label{eq2.54}
\end{equation}
$$\left.+\int\limits_0^xG(\lambda,x,t)[\theta s_2(\lambda,l)+s_2^*(\lambda,l)]f(t)dt\right\}.$$
Taking into account \eqref{eq2.52}, we have
$$G'(\lambda,l,t)s_2(\lambda,l)-G(\lambda,l,t)s'_2(\lambda,l)=s_2(\lambda,l)[s'_0(\lambda,l)s_2^*(\lambda,t)-s'_1(\lambda,l)s_1^*(\lambda,t)$$
$$+s'_2(\lambda,l)s_0^*(\lambda,t)]-s'_2(\lambda,l)[s_0(\lambda,l)s_2^*(\lambda,t)-s_1(\lambda,l)s_1^*(\lambda,t)+s_2(\lambda,l)s_0^*(\lambda,t)]$$
$$=s_2^*(\lambda,t)W_{2,0}(\lambda,l)-s_1^*(\lambda,t)W_{2,1}(\lambda,l)=-s_2^*(\lambda,t)s_1^*(\lambda,l)+s_1^*(\lambda,l)s_2^*(\lambda,l)$$
(see \eqref{eq2.32}), and analogously
$$G(\lambda,l,t)s'_1(\lambda,l)-G'(\lambda,l,t)s_1(\lambda,l)=-s_2^*(\lambda,t)s_0^*(\lambda,l)-s_0^*(\lambda,t)s_2^*(\lambda,l),$$
therefore
$$s_1(\lambda,x)[G'(\lambda,l,t)s_2(\lambda,l)-G(\lambda,l,t)s'_2(\lambda,l)]+s_2(\lambda,x)[G(\lambda,l,t)s'_1(\lambda,l)$$
$$-G'(\lambda,l,t)s_1(\lambda,l)]=s_2^*(\lambda,t)(-s_1(\lambda,x)s_1^*(\lambda,l)+s_2(\lambda,x)s_0^*(\lambda,l))$$
$$+s_2^*(\lambda,l)(s_1(\lambda,x)s_1^*(\lambda,t)-s_2(\lambda,x)s_0^*(\lambda,t))=s_2^*(\lambda,t)G(\lambda,x,l)-s_2^*(\lambda,l)G(\lambda,x,t).$$
Hence, according to \eqref{eq2.54}, we find
\begin{equation}
\begin{array}{ccc}
{\displaystyle Y(\lambda,x)=\frac i{\Delta_\theta(q,\lambda)}\left\{\int\limits_0^l[G(\lambda,x,l)s_2^*(\lambda,t)-G(\lambda,x,t)s_2^*(\lambda,l)-\theta G(\lambda,l,t)s_2(\lambda,x)]\right.}\\
{\displaystyle\left.\times f(t)dt+\int\limits_0^xG(\lambda,x,t)[\theta s_2(\lambda,l)+s_2^*(\lambda,l)]f(t)dt\right\}.}
\end{array}\label{eq2.55}
\end{equation}

\begin{lemma}\label{l2.9}
Resolvent $R_{L_q}(\lambda^3)=(L_q(\theta)-\lambda^3I)^{-1}$ of the operator $L_q$ \eqref{eq2.1} is
\begin{equation}
\begin{array}{ccc}
{\displaystyle(R_{L_q}(\lambda^3)f)(x)=\frac i{\Delta_\theta(q,\lambda)}\left\{\int\limits_0^l[G(\lambda,x,l)s_2^*(\lambda,t)-\theta G(\lambda,l,t)s_2(\lambda,x)]f(t)dt\right.}\\
{\displaystyle\left.+\theta\int\limits_0^xG(\lambda,x,t)s_2(\lambda,l)f(t)dt-\int\limits_x^lG(\lambda,x,t)s_2^*(\lambda,l)f(t)dt\right\}}
\end{array}\label{eq2.56}
\end{equation}
where $\Delta_\theta(q,\lambda)$ is the characteristic function \eqref{eq2.34}; $G(\lambda,x,t)$ is given by \eqref{eq2.52}, and the functions $\{s_p(\lambda,x)\}$, $\{s_p^*(\lambda,x)\}$ are given by the formulas \eqref{eq2.23}, \eqref{eq2.33}.
\end{lemma}
\vspace{5mm}

{\bf 2.5} Calculate the orthogonal projection $E_n$ onto the proper subspace corresponding to $\lambda_n(q,\theta)$

\begin{theorem}\label{t2.2}
Orthogonal projection $E_n$ onto the proper subspace corresponding to the eigenvalue $\lambda_n^3(q,\theta)$ of the operator $L_q(\theta)$, where $\lambda_n(q,\theta)$ is a real zero of the characteristic function $\Delta_\theta(q,\lambda)$ \eqref{eq2.34}, equals
\begin{equation}
(E_nf)(x)\stackrel{\rm def}{=}\lim\limits_{\lambda\rightarrow\lambda_n}(\lambda_n^3-\lambda^3)(R_{L_q}(\lambda^3)f)(x)=\int\limits_0^lf(t)\overline{\psi_n(q,\lambda,t)}dt\psi_n(q,\lambda,x),
\label{eq2.57}
\end{equation}
besides, $\psi_n(q,\lambda,x)$ is an eigenfunction \eqref{eq2.45} of the operator $L_q(\theta)$ and
\begin{equation}
a_n^2(\lambda)=-\frac{\theta s_2(\lambda,l)+s_2^*(\lambda,l)}{i(\lambda_n^3-\lambda^3)}s_2(\lambda,l).\label{eq2.58}
\end{equation}
\end{theorem}

P r o o f. Let $\lambda=\lambda_n(q,\theta)\in\mathbb{R}$ be a zero of $\Delta_\theta(q,\lambda)$ \eqref{eq2.34}, then \eqref{eq2.55} implies
\begin{equation}
(E_nf)(x)=b_n\int\limits_0^l[G(\lambda,x,l)s_2^*(\lambda,t)-G(\lambda,x,t)s_2^*(\lambda,l)-\theta G(\lambda,l,t)s_2(\lambda,x)]f(t)dt\label{eq2.59}
\end{equation}
where $b_n=-\lim\limits_{\lambda\rightarrow\lambda_n}i(\lambda_n^3-\lambda^3)/(\theta s_2(\lambda,l)+s_2^*(\lambda,l))$. Consider the function
\begin{equation}
F(\lambda,x,t)\stackrel{\rm def}{=}G(\lambda,x,l)s_2^*(\lambda,t)-G(\lambda,x,t)s_2^*(\lambda,l)-\theta G(\lambda,l,t)s_2(\lambda,x),\label{eq2.60}
\end{equation}
assuming that $\lambda$ is a zero of $\Delta_\theta(q,\lambda)$, i. e.,
\begin{equation}
\theta s_2(\lambda,l)+s_2^*(\lambda,l)=0.\label{eq2.61}
\end{equation}
Equation \eqref{eq2.52} implies
$$F(\lambda,x,0)=-s_2(\lambda,x)s_2^*(\lambda,l)-\theta s_2(\lambda,l)s_2(\lambda,x)=0,$$
due to \eqref{eq2.61}. Upon differentiating with respect to $t$ ($\partial_t=d/dt$), we obtain
$$\partial_tF(\lambda,x,t)=G(\lambda,x,l)\partial_ts_2^*(\lambda,t)-\partial_tG(\lambda,x,t)s_2^*(\lambda,l)-\theta\partial_tG(\lambda,l,t)\cdot s_2(\lambda,x),$$
and using \eqref{eq2.52} and \eqref{eq2.61} we have
$$\left.\partial_tF(\lambda,x,t)\right|_{t=0}=s_1(\lambda,x)s_2^*(\lambda,l)+\theta s_1(\lambda,l)s_2(\lambda,x)=\theta u(\lambda,x)$$
where
\begin{equation}
u(\lambda,x)\stackrel{\rm def}{=}s_2(\lambda,x)s_1(\lambda,l)-s_2(\lambda,l)s_1(\lambda,x).\label{eq2.62}
\end{equation}
Upon differentiating again with respect to $t$, we find that
$$\partial_t^2F(\lambda,x,t)=G(\lambda,x,l)\partial_t^2s_2^*(\lambda,t)-\partial_t^2G(\lambda,x,t)s_2^*(\lambda,l)-\theta\partial_t^2G(\lambda,l,t)\cdot s_2(\lambda,x),$$
and thus
$$\left.\partial_t^2F(\lambda,x,t)\right|_{t=0}=G(\lambda,x,l)-s_0(\lambda,x)s_2^*(\lambda,l)-\theta s_0(\lambda,l)s_2(\lambda,x)$$
$$=s_2(\lambda,x)[s_0^*(\lambda,l)-\theta s_0(\lambda,l)]-s_1(\lambda,x)s_1^*(\lambda,l).$$
Notice that
\begin{equation}
\frac{s_1^*(\lambda,l)}{s_2(\lambda,l)}=\frac{s_0^*(\lambda,l)-\theta s_0(\lambda,l)}{s_1(\lambda,l)}(\stackrel{\rm def}{=}R).\label{eq2.63}
\end{equation}
This equality, due to \eqref{eq2.61} and \eqref{eq2.52}, is equivalent to the relation
$$s_1^*(\lambda,l)s_1(\lambda,l)=s_2(\lambda,l)s_0^*(\lambda,l)+s_0(\lambda,l)s_2^*(\lambda,l)\Leftrightarrow G(\lambda,l,l)=0.$$
Using \eqref{eq2.63}, we obtain that
$$\left.\partial_t^2F(\lambda,x,t)\right|_{t=0}=Ru(\lambda,x)$$
where $u(\lambda,x)$ is given by \eqref{eq2.62}. So, for all $x\in[0,l]$, function $F(\lambda,x,t)$ \eqref{eq2.60} is a solution to the Cauchy problem
$$i\partial_t^3F(\lambda,x,t)=(q(t)-\lambda^3)F(\lambda,x,t);$$
$$\left.F(\lambda,x,t)\right|_{t=0}=0;\quad\left.\partial_tF(\lambda,x,t)\right|_{t=0}=\theta u(\lambda,x);\quad\left.\partial_t^2F(\lambda,x,t)\right|_{t=0}=Ru(\lambda,x);$$
therefore
$$F(\lambda,x,t)=\theta u(\lambda,x)s_1^*(\lambda,t)+Ru(\lambda,x)s_2^*(\lambda,t)=u(\lambda,x)[\theta s_1^*(\lambda,t)+Rs_2^*(\lambda,t)],$$
and according to \eqref{eq2.61}, \eqref{eq2.63},
$$F(\lambda,x,t)=\frac{u(\lambda,x)}{s_2(\lambda,l)}[s_2^*(\lambda,t)s_1^*(\lambda,l)-s_1^*(\lambda,t)s_2^*(\lambda,l)]=\frac{u(\lambda,x)\overline{u(\lambda,t)}}{s_2(\lambda,l)}$$
($\forall\lambda\in\mathbb{R}$). And since $u(\lambda,x)=a_n\psi_n(q,\lambda,x)$ (see \eqref{eq2.45}) where $\lambda=\lambda_n(q,\theta)$, then
$$F(\lambda,x,t)=\frac{a_n^2(\lambda)}{s_2(\lambda,l)}\overline{\psi_n(q,\lambda,t)}\psi_n(q,\lambda,x),$$
and thus equality \eqref{eq2.59} becomes
$$(E_nf)(x)=\frac{b_na_n^2}{s_2(\lambda,l)}\langle f(t),\psi_n(q,\lambda,t)\rangle\psi_n(q,\lambda,x).$$
Since $E_n$ is an orthogonal projection, then $b_na_n^2=s_2(x,l)$, which proves \eqref{eq2.58}. $\blacksquare$

\begin{remark}\label{r2.2}
Zeros $\lambda=\lambda_n(q,\theta)$ of the characteristic function $\Delta_\theta(q,\lambda)$ \eqref{eq2.34} are the only singularities of the resolvent $R_{L_q}(\lambda^3)$ \eqref{eq2.56}. Orthogonal projections $E_n$ onto the proper subspaces corresponding to $\lambda_n(q,\theta)$ are one-dimensional \eqref{eq2.56}. And since ${\displaystyle\sum\limits_nE_n=I}$, the functions $\{\psi_n(q,\lambda,x)\}$ \eqref{eq2.45} form the orthonormal basis of the space $L^2(0,l)$.
\end{remark}

\section{The main equations}\label{s3}

{\bf 3.1} By $w(\lambda,x)$, we denote  the solution to equation \eqref{eq2.2},
\begin{equation}
w(\lambda,x)\stackrel{\rm def}{=}\widehat{s}_2(\lambda,x)\widehat{s}_1(\lambda,l)-\widehat{s}_1(\lambda,x)\widehat{s}_2(\lambda,l)\quad(\widehat{s}_p(\lambda,x)\stackrel{\rm def}{=}(i\lambda)^ps_p(\lambda,x))\label{eq3.1}
\end{equation}
which satisfies the boundary conditions $y(0)=0$ and $y(l)=0$. And let $e_k(\lambda,x)$ be obtained from $e^{i\lambda\zeta_kx}$ via the transformation operator \eqref{eq2.10},
\begin{equation}
e_k(\lambda,x)\stackrel{\rm def}{=}(I+T_\lambda)e^{i\lambda\zeta_kx}=e^{i\lambda\zeta_kx}+\int\limits_0^xT(\lambda,x,t)e^{i\lambda\zeta_kt}dt\quad(1\leq k\leq3)\label{eq3.2}
\end{equation}
where the kernel $T(\lambda,x,t)$ is given by \eqref{eq2.19}. Then $w(\lambda,x)$ \eqref{eq3.1} is written as
\begin{equation}
w(\lambda,x)=B_1(\lambda)e_1(\lambda,x)+B_2(\lambda)e_2(\lambda,x)+B_3(\lambda)e_3(\lambda,x)\label{eq3.3}
\end{equation}
where
\begin{equation}
\begin{array}{ccc}
{\displaystyle B_1(\lambda)\stackrel{\rm def}{=}\frac13(\widehat{s}_1(\lambda,l),\widehat{s}_2(\lambda,l));\quad B_2(\lambda)\stackrel{\rm def}{=}\frac13(\zeta_2\widehat{s}_1(\lambda,l)-\zeta_3\widehat{s}_2(\lambda,l));}\\
{\displaystyle B_3(\lambda)\stackrel{\rm def}{=}\frac13(\zeta_3\widehat{s}_1(\lambda,l)-\zeta_2\widehat{s}_2(\lambda,l)),}
\end{array}\label{eq3.4}
\end{equation}
besides, $\widehat{B}_k(\lambda\zeta_2)=\widehat{B}_{k'}(\lambda)$ ($k'=(k+1)\mod3$). Define the functions
\begin{equation}
\omega_p(\lambda,x)\stackrel{\rm def}{=}w(\lambda,x)/B_p(\lambda)\quad(1\leq p\leq3),\label{eq3.5}
\end{equation}
then  \eqref{eq3.3} implies
\begin{equation}
\omega_1(\lambda,x)=e_1(\lambda,x)+c_2(\lambda)e_2(\lambda,x)+c_3(\lambda)e_3(\lambda,x)\label{eq3.6}
\end{equation}
where
\begin{equation}
c_2(\lambda)\stackrel{\rm def}{=}B_2(\lambda)/B_1(\lambda);\quad c_3(\lambda)\stackrel{\rm def}{=}B_3(\lambda)/B_1(\lambda).\label{eq3.7}
\end{equation}

\begin{lemma}\label{l3.1}
The functions $\{c_p(\lambda)\}$ \eqref{eq3.7} satisfy the following relations:
\begin{equation}
\begin{array}{lll}
({\rm i})&c_2(\lambda)c_2(\lambda\zeta_2)c_2(\lambda\zeta_3)=1;\\
({\rm ii})&c_2(\lambda\zeta_3)c_3(\lambda)=1.
\end{array}\label{eq3.8}
\end{equation}
\end{lemma}

Proof of \eqref{eq3.8} follows from the equalities
$$c_2(\lambda\zeta_2)=B_3(\lambda)/B_2(\lambda);\quad c_2(\lambda\zeta_3)=B_1(\lambda)/B_3(\lambda).$$
\vspace{5mm}

{\bf 3.2} For the Wronskians $\{W_{k,s}^e(\lambda,x)\}$ of the functions $\{e_k(\lambda,x)\}$ \eqref{eq3.2}
\begin{equation}
W_{k,s}^e(\lambda,x)\stackrel{\rm def}{=}e_k(\lambda,x)e'_s(\lambda,x)-e_s(\lambda,x)e'_k(\lambda,x)\quad(1\leq k,s\leq3),\label{eq3.9}
\end{equation}
the following statement holds.

\begin{lemma}\label{l3.2}
Wronskians $\{W_{k,s}^e(\lambda,x)\}$ \eqref{eq3.9} have the following representation:
\begin{equation}
W_{1,2}^e(\lambda,x)=\sqrt3\lambda\zeta_3e_2^*(\lambda,x);\quad W_{2,3}^e(\lambda,x)=\sqrt3\lambda\zeta_1e_1^*(\lambda);\quad W_{3,1}^e(\lambda,x)=\sqrt3\zeta_2e_3^*(\lambda,x)\label{eq3.10}
\end{equation}
\end{lemma}

P r o o f. The function
$$W_{1,2}^e(\lambda,x)=e_1(\lambda,x)e'_2(\lambda,x)-e_2(\lambda,x)e'_1(\lambda,x)\quad(\left.W_{1,2}^e(\lambda,x)\right|_{x=0}=i\lambda(\zeta_2-\zeta_1))$$
upon differentiating becomes
$$(W_{1,2}^e(\lambda,x))'=e_1(\lambda,x)e''_2(\lambda,x)-e_2(\lambda,x)e''_1(\lambda,x)\quad(\left.(W_{1,2}^e(\lambda,x))'\right|_{x=0}=(i\lambda)^2(\zeta_3-\zeta_1)),$$
therefore
$$(W_{1,2}^e(\lambda,x))''=e'_1(\lambda,x)e''_2(\lambda,x)-e'_2(\lambda,x)e''_1(\lambda,x)\quad(\left.(W_{1,2}^e(\lambda,x))''\right|_{x=0}=(i\lambda)^3(\zeta_3-\zeta_2)).$$
Hence it follows that $W_{1,2}^e(\lambda,x)$ is the solution to the equation
$$iy'''(x)-q(x)y(x)=-\lambda^3y(x)$$
satisfying the initial data $y(0)=\sqrt3\lambda\zeta_3$, $y'(0)=-i\sqrt3\lambda^2\zeta_2$; $y''(0)=-\sqrt3\lambda^3\zeta_1$ (since $\zeta_2-\zeta_1=-i\sqrt3\zeta_3$, $\zeta_3-\zeta_1=i\sqrt3\zeta_2$, $\zeta_3-\zeta_2=-i\sqrt3\zeta_1$), and thus $W_{1,2}^e(\lambda,x)=\sqrt3\lambda\zeta_3e_2^*(\lambda,x)$. Other relations in \eqref{eq3.10} follow from the proven one after the substitutions $\lambda\rightarrow\lambda\zeta_2$, $\lambda\rightarrow\lambda\zeta_3$. $\blacksquare$

Boundary conditions $y(0)=0$ and $y'(0)=\overline{\theta}y'(l)$ \eqref{eq1.18} for $\omega_1(\lambda,x)$ \eqref{eq3.6} yield the following equation system:
\begin{equation}
\left\{
\begin{array}{lll}
e_1(\lambda,0)+c_2(\lambda)e_2(\lambda,0)+c_3(\lambda)e_3(\lambda,0)=0;\\
e'_1(\lambda,0)+c_2(\lambda)e'_2(\lambda,0)+c_3(\lambda)e'_3(\lambda,0)=i\lambda\overline{\theta}\widehat{s}_2^*(\lambda,l)/B_1(\lambda)
\end{array}\right.\label{eq3.11}
\end{equation}
since $\omega'_1(\lambda,x)=(i\lambda)^3W_{1,2}(\lambda,l)/B_1(\lambda)$ and $W_{1,2}(\lambda,l)=s_2^*(\lambda,l)$ \eqref{eq2.32}. Upon multiplying the first equality in \eqref{eq3.11} by $e'_p(\lambda,0)$ and the second, by $e_p(\lambda,0)$, and subtracting, we have
$$\left\{
\begin{array}{lll}
c_2(\lambda)W_{2,1}^e(\lambda,0)+c_3(\lambda)W_{3,1}^e(\lambda,0)=-n(\lambda)e_1(\lambda,0);\\
W_{1,2}^e(\lambda,0)+c_3(\lambda)W_{3,2}^e(\lambda,0)=-n(\lambda)e_0(\lambda,0);\\
W_{1,3}^e(\lambda,0)+c_2(\lambda)W_{2,3}^e(\lambda,0)=-n(\lambda)e_3(\lambda,0);
\end{array}\right.$$
where $n(\lambda)\stackrel{\rm def}{=}i\lambda\theta\widehat{s}_2^*(\lambda,l)/B_1(\lambda)$. Using \eqref{eq3.10}, we obtain
$$\left\{
\begin{array}{lll}
{\displaystyle-s_3c_2(\lambda)e_2^*(\lambda,0)+\zeta_2c_3(\lambda)e_3^*(\lambda,0)=-\frac{n(\lambda)}{\sqrt3\lambda}e_1(\lambda,0);}\\
{\displaystyle\zeta_3e_2^*(\lambda,0)-\zeta_1c_3(\lambda)e_1^*(\lambda,0)=-\frac{n(\lambda)}{\sqrt3\lambda}e_2(\lambda,0);}\\
{\displaystyle-\zeta_2e_3^*(\lambda,0)+\zeta_1c_2(\lambda)e_1^*(\lambda,0)=-\frac{n(\lambda)}{\sqrt3\lambda}e_3(\lambda,0).}
\end{array}\right.$$
Rewrite this system in the matrix form:
\begin{equation}
T(\lambda)E^*(\lambda,0)=-\frac{n(\lambda)}{\sqrt3\lambda}E(\lambda,l)\label{eq3.12}
\end{equation}
where
$$E(\lambda,0)\stackrel{\rm def}{=}\left[
\begin{array}{ccc}
e_1(\lambda,0)\\
e_2(\lambda,0)\\
e_3(\lambda,0)
\end{array}\right];\quad T(\lambda)\stackrel{\rm def}{=}\left[
\begin{array}{ccc}
0&-\zeta_3c_2(\lambda)&\zeta_2c_3(\lambda)\\
-\zeta_2c_3(\lambda)&\zeta_3&0\\
\zeta_1c_2(\lambda)&0&-\zeta_2
\end{array}
\right]$$
where $E^*(\lambda,0)=\overline{E(\overline{\lambda},0)}$ (see \eqref{eq2.31}. Applying the ``$*$''-operation \eqref{eq2.31} to equality \eqref{eq3.12}, we have
$$T^+(\lambda)E(\lambda,0)=-\frac{n^*(\lambda)}{\sqrt3\lambda}E^*(\lambda,0)$$
where $T^+(\lambda)$ is obtained from the matrix $T(\lambda)$ via application of ``$*$'' to each its element. Hence and from \eqref{eq3.12} it follows that
\begin{equation}
T(\lambda)T^+(\lambda)E(\lambda,0)=\frac{n^*(\lambda)n(\lambda)}{3\lambda^2}E(\lambda,0).\label{eq3.13}
\end{equation}
Upon equating the first elements of the vectors in this equality, we obtain
$$[\zeta_3c_1(\lambda)c_3^*(\lambda)+\zeta_2c_3(\lambda)c_2^*(\lambda)]e_1(\lambda,0)-c_2(\lambda)e_2(\lambda,0)-c_3(\lambda)e_3(\lambda,0)$$
$$=\frac{n(\lambda)n^*(\lambda)}{3\lambda^2}e_1(\lambda,0).$$
And taking into account the first relation in \eqref{eq3.11}, we have
\begin{equation}
\zeta_3c_2(\lambda)c_3^*(\lambda)+\zeta_2c_3(\lambda)c_2^*(\lambda)+1=\frac{\widehat{s}_2(\lambda,l)\widehat{s}_2^*(\lambda,l)}{3B_1(\lambda)B_1^*(\lambda)}.\label{eq3.14}
\end{equation}
Equating the second and the third components of the vectors in \eqref{eq3.13} again leads to equality \eqref{eq3.14}.

\begin{lemma}\label{l3.3}
For the functions $c_2(\lambda)$, $c_3(\lambda)$ \eqref{eq3.7}, relation \eqref{eq3.14} holds where $\widehat{s}_2(\lambda,l)=(i\lambda)^2s_2(\lambda,l)$ and $B_1(\lambda)$ is from \eqref{eq3.4}.
\end{lemma}

\begin{remark}\label{r3.1}
Equalities \eqref{eq3.8} follow from invariancy of equation \eqref{eq2.8} relative to the transform $\lambda\rightarrow\lambda\zeta_2$, and relation \eqref{eq3.14} is the conservation law (`unitarity' property) that follow from self-adjoint boundary conditions.
\end{remark}
\vspace{5mm}

{\bf 3.3} Relation \eqref{eq3.6} implies that
\begin{equation}
\left\{
\begin{array}{lll}
\omega_1(\lambda,x)=e_1(\lambda,x)+c_2(\lambda)e_2(\lambda,x)+c_3(\lambda)e_3(\lambda,x);\\
\omega'_1(\lambda,x)=e'_1(\lambda,x)+c_2(\lambda)e'_2(\lambda,x)+c_3(\lambda)e'_3(\lambda,x).
\end{array}\right.\label{eq3.15}
\end{equation}
Upon multiplying the first of the equalities in \eqref{eq3.15} by $e'_3(\lambda,x)$ and the second, by $e_3(\lambda,x)$ and subtracting, we obtain
$$\omega_{1,3}(\lambda,x)=W_{1,3}^e(\lambda,x)+c_2(\lambda)W_{2,3}^e(\lambda,x)$$
where $\omega_{1,3}(\lambda,x)$ is the Wronskian of $\omega_1(\lambda,x)$ and $e_3(\lambda,x)$,
\begin{equation}
\omega_{p,s}(\lambda,x)\stackrel{\rm def}{=}\omega_p(\lambda,x)e'_s(\lambda,x)-\omega'_p(\lambda,x)e_s(\lambda,x)\quad(1\leq p,s\leq3).\label{eq3.16}
\end{equation}
Hence, due to \eqref{eq3.10}, follows the statement.

\begin{lemma}\label{l3.4} For all $\lambda\in\mathbb{C}$, the following equalities hold:
\begin{equation}
\begin{array}{lll}
({\rm i})&{\displaystyle c_2(\lambda)e_1^*(\lambda,x)=\zeta_2e_3^*(\lambda,x)+\frac{\zeta_1}{\sqrt3\lambda}\omega_{1,3}(\lambda,x);}\\
({\rm ii})&{\displaystyle c_2(\lambda\zeta_2)e_3^*(\lambda,x)=\zeta_2e_2^*(\lambda,x)+\frac{\zeta_3}{\sqrt3\lambda}\omega_{2,1}(\lambda,x);}\\
({\rm iii})&{\displaystyle c_2(\lambda\zeta_3)e_2^*(\lambda,x)=\zeta_2e_1^*(\lambda,x)+\frac{\zeta_2}{\sqrt3\lambda}\omega_{3,2}(\lambda,x);}
\end{array}\label{eq3.17}
\end{equation}
where $\{\omega_{p,s}(\lambda,x)\}$ are given by \eqref{eq3.16}.
\end{lemma}

Similarly, multiplying the first relation in \eqref{eq3.15} by $e'_2(\lambda,x)$ and the second, by $e_2(\lambda,x)$, and subtracting, we arrive at the equality
$$\omega_{1,2}(\lambda,x)=W_{1,2}^e(\lambda,x)+c_3(\lambda)W_{3,2}^e(\lambda,x),$$
which gives the following statement.

\begin{lemma}\label{l3.5}
For all $\lambda\in\mathbb{C}$, the following relations take place:
\begin{equation}
\begin{array}{lll}
({\rm i})&{\displaystyle c_3(\lambda)e_1^*(\lambda)=\zeta_3e_2^*(\lambda,x)-\frac{\zeta_1}{\sqrt3\lambda}\omega_{1,2}(\lambda,x);}\\
({\rm ii})&{\displaystyle c_3(\lambda\zeta_2)e_3^*(\lambda)=\zeta_3e_1^*(\lambda,x)-\frac{\zeta_3}{\sqrt3\lambda}\omega_{2,3}(\lambda,x);}\\
({\rm iii})&{\displaystyle c_3(\lambda\zeta_3)e_2^*(\lambda,x)=\zeta_3e_3^*(\lambda,x)-\frac{\zeta_2}{\sqrt3\lambda}\omega_{3,1}(\lambda,x).}
\end{array}\label{eq3.18}
\end{equation}
\end{lemma}

\begin{remark}\label{r3.2}
Equalities \eqref{eq3.17} and \eqref{eq3.18} are equivalent and all of them could be derived from any of the relations \eqref{eq3.17} (or \eqref{eq3.18}). Show that {\rm (i)} \eqref{eq3.17} and {\rm(ii)} \eqref{eq3.18} coincide. Rewrite {\rm(i)} \eqref{eq3.17} as
$$B_2(\lambda)e_1^*(\lambda,x)-\zeta_2B_1(\lambda)e_3^*(\lambda,x)=\frac{\zeta_1}{\sqrt3\lambda}[\omega(\lambda,x)e'_3(\lambda,x)-\omega'(\lambda,x)e_3(\lambda,x)],$$
and equality {\rm(ii)} \eqref{eq3.18}, in view of $c_3(\lambda\zeta_2)=B_1(\lambda)/B_2(\lambda)$, correspondingly, as
$$B_1(\lambda)e_3^*(\lambda)-\zeta_2B_2(\lambda)e_1^*(\lambda)=-\frac{\zeta_3}{\sqrt3\lambda}[\omega(\lambda,x)e'_3(\lambda,x)-\omega'(\lambda,x)e_3(\lambda,x)],$$
hence it follows the identity of the equalities {\rm(i)} \eqref{eq3.17} and {\rm(ii)} \eqref{eq3.18}.
\end{remark}
\vspace{5mm}

{\bf 3.4} Relations \eqref{eq2.20} and \eqref{eq3.2} imply that $\{e_k(\lambda,x)\}$ are solutions to the integral equations
\begin{equation}
e_k(\lambda,x)=e^{i\lambda\zeta_kx}+i\int\limits_0^xK_1(\lambda,x,t)q(t)e_k(\lambda,t)dt\quad(1\leq k\leq3)\label{eq3.19}
\end{equation}
($K_1(\lambda,x,t)$ is given by \eqref{eq2.6}), and for the functions
\begin{equation}
E_k(\lambda,x)\stackrel{\rm def}{=}e_k(\lambda,x)e^{-i\lambda\zeta_kx}\quad(1\leq k\leq3)\label{eq3.20}
\end{equation}
the following relations hold:
\begin{equation}
E_k(\lambda,x)=1+i\int\limits_0^xe^{-i\lambda\zeta_k(x-t)}K_1(\lambda,x,t)q(t)E_k(\lambda,t)dt\quad(1\leq k\leq3).\label{eq3.21}
\end{equation}
Since
$$K_1(\lambda,x,t)e^{-i\lambda\zeta_1(x-t)}=\frac1{3(i\lambda)^2}\{1+\zeta_2e^{i\lambda(\zeta_2-\zeta_1)(x-t)}+\zeta_3e^{i\lambda(\zeta_3-\zeta_1)(x-t)}\}$$
and
$$i\lambda(\zeta_2-\zeta_1)=\frac12(3\beta-\alpha\sqrt3)+\frac i2(-3\alpha-\sqrt3\beta);$$
$$ i\lambda(\zeta_3-\zeta_1)=\frac12(3\beta+\alpha\sqrt3)+\frac i2(-3\alpha+\sqrt3\beta)$$
($\lambda=\alpha+i\beta$; $\alpha$, $\beta\in\mathbb{R}$), then every exponent in the domain $\sqrt3\beta-\alpha<0$, $\sqrt3\beta+\alpha<0$ (coinciding with $\Omega_1$ \eqref{eq1.14}) is less then $1$ in modulo and vanishes when $\lambda\rightarrow\infty$ ($\lambda\in\Omega_1$).

\begin{lemma}\label{l3.6}
For the functions $\{e_k(\lambda,x)\}$ \eqref{eq3.2} in the sectors $\{\Omega_k\}$ \eqref{eq1.14}, the following formulas are true:
\begin{equation}
\begin{array}{lll}
{\rm(a)}&{\displaystyle E_k(\lambda,x)=1-\frac{i\zeta_k}{3\lambda^2}\int\limits_0^xq(t)dt+O\left(\frac1{\lambda^4}\right);}\\
{\rm(b)}&{\displaystyle e^{-i\lambda\zeta_kx}e'_k(\lambda,x)=i\lambda\zeta_k+\frac1{3\lambda\zeta_k}\int\limits_0^xq(t)dt+O\left(\frac1{\lambda^3}\right)}
\end{array}\quad(\lambda\in\Omega_k;|\lambda|\gg1)\label{eq3.22}
\end{equation}
where $\{E_k(\lambda,x)\}$ are given by \eqref{eq3.20}.
\end{lemma}

Analogously, equation \eqref{eq1.16} and
$$\widehat{s}_p(\lambda,x)=s_p(i\lambda x)+i\int\limits_0^xK_1(\lambda,x,t)q(t)\widehat{s}_p(\lambda,t)dt\quad(0\leq p\leq2)$$
yield the statement.

\begin{lemma}\label{l3.7}
Inside the sectors $\{\Omega_k\}$\eqref{eq1.14} for $\{\widehat{s}_p(\lambda,x)\}$, the following relations hold:
\begin{equation}
\begin{array}{lll}
({\rm a})&{\displaystyle e^{-i\lambda\zeta_kx}3\zeta_k^p\widehat{s}_p(\lambda,x)=1-\frac{i\zeta_k}{3\lambda^2}\int\limits_0^xq(t)dt+O\left(\frac1{\lambda^2}\right);}\\
({\rm b})&{\displaystyle e^{-i\lambda\zeta_kx}3\zeta_k^p\widehat{s}'_p(\lambda,x)=i\lambda\zeta_k+\frac{\zeta_k^2}{3\lambda}\int\limits_0^xq(t)dt+O\left(\frac1{\lambda^3}\right)}
\end{array}\quad(\lambda\in\Omega_k,|\lambda|\gg1)\label{eq3.23}
\end{equation}
($1\leq k\leq3$, $0\leq p\leq2$).
\end{lemma}
\vspace{5mm}

{\bf 3.5} Equalities \eqref{eq3.17} and \eqref{eq3.18} generate problems of a jump in adjacent sectors of half-plane. Let $\{\Omega_k^-\}$ be the sectors
\begin{equation}
\Omega_k^-\stackrel{\rm def}{=}\{\lambda\in\mathbb{C}:-\lambda\in\Omega_k\}\quad(1\leq k\leq3)\label{eq3.24}
\end{equation}
obtained from $\{\Omega_k\}$ \eqref{eq1.14} after the substitution $\lambda\rightarrow-\lambda$. Rewrite relation (i) \eqref{eq3.18} in terms of $\{E_k(\lambda,x)\}$ \eqref{eq3.20},
\begin{equation}
\zeta_3e^{i\lambda(\zeta_1-\zeta_2)x}c_3^*(\lambda)E_1(\lambda,x)=E_2(\lambda,x)-f_{1,2}(\lambda,x)\quad\left(f_{1,2}(\lambda,x)\stackrel{\rm def}{=}\frac{\zeta_3}{\sqrt3\lambda}e^{-i\lambda\zeta_2x}\omega_{1,2}^k(\lambda,x)\right).\label{eq3.25}
\end{equation}
Entire functions $\{E_k(\lambda)\}$ tend to 1 as $\lambda\rightarrow\infty$ ($\lambda\in\Omega_k$), due to (a) \eqref{eq3.2}. Show that $f_{1,2}(\lambda,x)$ \eqref{eq3.25} also tends to 1 when $\lambda\rightarrow\infty$ and $\lambda\in\Omega_1\cap\Omega_3^-$.

\begin{picture}(200,200)
\put(0,100){\vector(1,0){200}}
\put(100,0){\vector(0,1){200}}
\put(100,100){\vector(3,-1){100}}
\put(100,100){\vector(-3,-1){100}}
\qbezier(100,130)(130,110)(135,90)
\qbezier(140,85)(100,70)(60,85)
\qbezier(63,90)(80,120)(100,128)
\qbezier(40,80)(70,50)(100,60)
\qbezier(45,83)(75,53)(100,63)
\put(110,190){$il_{\zeta_1}$}
\put(110,140){$E_3(\lambda,x)$}
\put(50,140){$E_2(\lambda,x)$}
\put(103,102){$\Omega_3$}
\put(180,110){$l_{\zeta_1}$}
\put(180,62){$il_{\zeta_3}$}
\put(83,85){$\Omega_1$}
\put(6,77){$il_{\zeta_2}$}
\put(83,102){$\Omega_2$}
\put(110,50){$E_1(\lambda,x)$}
\put(30,50){$\Omega_1\cap\Omega'_3$}
\put(35,25){$f_{1,2}(\lambda,x)$}
\end{picture}

\hspace{20mm} Fig. 2.

Since
$$f_{1,2}(\lambda,x)=\frac{\zeta_3e^{-i\lambda\zeta_2x}}{\sqrt3\lambda B_1^*(\lambda)}\{w^*(\lambda,x),e_2^*(\lambda,x)\}$$
$$=\frac{\zeta_3}{\sqrt3\lambda B_1^*(\lambda)}\{e^{i\lambda\zeta_1x}w^*(\lambda,x)\cdot e^{i\lambda\zeta_3x}e_2^{*'}(\lambda,x)$$
$$-e^{i\lambda\zeta_1x}w^{*'}(\lambda,x)\cdot e^{i\lambda\zeta_3x}e_2^*(\lambda,x)\},$$
then, taking into account \eqref{eq3.1} and asymptotic \eqref{eq3.22}, \eqref{eq3.23} (which is true for $\lambda\in\Omega_1\cap\Omega_3^-$), we obtain
$$f_{1,2}(\lambda,x)=\frac{\zeta_3}{\sqrt3\lambda B_1^*(\lambda)}\left\{\frac13(s_1^*(\lambda,l)-s_2^*(\lambda,l))\left(1+O\left(\frac1{\lambda^2}\right)\right)\left(-i\lambda\zeta_3+O\left(\frac1\lambda\right)\right)\right.$$
$$\left.+\frac13(
s_1^*(\lambda,l)-s_2^*(\lambda,l))\left(i\lambda+O\left(\frac1\lambda\right)\right)\left(1+O\left(\frac1{\lambda^2}\right)\right)\right\}=\frac{\zeta_3}{3\sqrt3\lambda B_1^*(\lambda)}(s_1^*(\lambda,l)$$
$$-s_2^*(\lambda,l))\left(i\lambda(\zeta_1-\zeta_3)+O\left(\frac1\lambda\right)\right).$$
Using \eqref{eq3.4} and $\zeta_1-\zeta_2=-i\sqrt3\zeta_2$, we arrive at the equality ${\displaystyle f_{1,2}(\lambda,x)=1+O\left(\frac1\lambda\right)}$ ($|\lambda|\gg1$, $\lambda\in\Omega_1\cap\Omega_3^-$).

Function $f_{1,2}(\lambda,x)$ \eqref{eq3.23} can have singularities at zeros of the function $B_1^*(\lambda)$. Relation $B_1^*(\lambda)=0\Leftrightarrow B_1(\overline{\lambda})=0$ implies that the set of zeros of $B_1^*(\lambda)$ is obtained after complex conjugation of zeros of $B_1(\lambda)$ \eqref{eq3.4}. Using $\widehat{s}_k(\lambda,x)=(I-iK_\lambda)s_k(i,\lambda x)$ \eqref{eq2.18} ($K_\lambda$ are given by \eqref{eq2.7}, we obtain the equation for $B_1(\lambda,x)\stackrel{\rm def}{=}\widehat{s}_1(\lambda,x)-\widehat{s}_2(\lambda,x)$:
$$B_1(\lambda,x)=s_1(i\lambda x)-s_2(i\lambda x)-i\int\limits_0^xK_1(\lambda,x,t)q(t)B_1(\lambda,t)dt$$
which, due to \eqref{eq2.6} and
$$s_1(i\lambda x)-s_2(i\lambda x)=\frac13(\zeta_3-\zeta_2)(e^{i\lambda\zeta_2x}-e^{i\lambda\zeta_3x})=\frac{2i}{\sqrt3}e^{-\frac{i\lambda x}2}\sh\frac{\sqrt3\lambda x}2,$$
we can rewrite as
\begin{equation}
B_1(\lambda,x)=\frac{2i}{\sqrt3}e^{-\frac{i\lambda x}2}\sh\frac{\sqrt3\lambda x}2+\frac i{\lambda^3}\int\limits_0^xs_2(i\lambda(x-t))q(t)B_1(\lambda,t)dt.\label{eq3.26}\
\end{equation}
For $q\equiv0$, the equality $B_1(\lambda,l)=0$ is equivalent to ${\displaystyle\sh\frac{\sqrt3\lambda l}2=0}$, and thus $\lambda=2\pi in/\sqrt3l$. So, the set of zeros of $B_1(\lambda,l)$, for $q\equiv0$, is
\begin{equation}
\Lambda_0\stackrel{\rm def}{=}\left\{\mu_n(0)=\frac{2\pi in}{\sqrt3 l}:n\in\mathbb{Z}\right\},\label{eq3.27}
\end{equation}
besides, $\overline{\mu_n(0)}=-\mu_n(0)=\mu_{-n}(0)$.

Equation \eqref{eq3.26} implies that the solution $B_1(\lambda,x)$ has the property $\overline{B_1(\lambda,x)}=B_1(-\overline{\lambda},x)$, hence it follows that zeros of $B_1(\lambda,x)$ are included in pairs $(\lambda,-\overline{\lambda})$ of points symmetric relative to the imaginary axis $i\mathbb{R}$.

Consider the circles $\mathbb{T}_R\stackrel{\rm def}{=}\{\lambda\in\mathbb{C}:|\lambda|=R\}$ of the radius $R$ ($0<R$, $R\gg1$), assuming that $\mathbb{T}_R\cap\Lambda_0=\emptyset$ ($\Lambda_0$ is given by \eqref{eq3.27}). By $\mathbb{T}_R^k\stackrel{\rm def}{=}\mathbb{T}_R\cap\Omega_k$, we denote the arc of the circle $\mathbb{T}_R$ lying in the sector $\Omega_k$ \eqref{eq1.14} ($1\leq k\leq3$). Show that modulo of the first summand in \eqref{eq2.26} is greater than modulo of the second summand on each arc $\mathbb{T}_R^k$ when $R\gg1$. Multiply equation \eqref{eq3.26} by $e^{-i\lambda\zeta_2x}$, then for $A_1(\lambda,x)=B_1(\lambda,x)e^{-i\lambda x}$ we obtain
\begin{equation}
A_1(\lambda,x)=\frac i{\sqrt3}\left(1-e^{i\lambda(\zeta_3-\zeta_2)x}\right)+\frac i{\lambda^2}\int\limits_0^xe^{-i\lambda\zeta_2(x-t)}s_2(i\lambda(x-t))q(t)A_1(\lambda,t)dt.\label{eq3.28}
\end{equation}
Since $i\lambda(\zeta_3-\zeta_2)x=\alpha\sqrt3+i\sqrt3\beta$ ($\lambda=\alpha+i\beta$; $\alpha$, $\beta\in\mathbb{R}$), then, for $\alpha<0$,
$$\frac i{\sqrt3}\left(1-e^{i\lambda(\zeta_3-\zeta_2)x}\right)>\frac1{\sqrt3}\left(1-e^{\alpha\sqrt3}\right)\quad(\lambda\in\Omega_2).$$
In the sector $\Omega_2$, ${\displaystyle\left|\frac1{\lambda^2}e^{-i\lambda\zeta_2(x-t)}s_2(i\lambda(x-t))\right|<\frac1{|\lambda|^2}}$, therefore
$$\left|\frac i{\lambda^2}\int\limits_0^xe^{-i\lambda\zeta_2(x-t)}s_2(i\lambda(x-t))q(t)A_1(\lambda,t)dt\right|\leq\frac{1+e^{\alpha\sqrt3}}{\sqrt3|\lambda|^2}\sigma\exp\left(\frac\sigma{|
\lambda|^2}\right)\quad|\lambda|\gg1$$
where ${\displaystyle\sigma=\int\limits_0^l|q(t)|dt}$ (see \eqref{eq2.17}) and consequently, for $|\lambda|\gg1$ ($\lambda\in\Omega_2$), the right-hand side of this inequality does not exceed $M/|\lambda|^2$ ($M>0$). Thus, on the arc $\mathbb{T}_R^2$, the following inequality holds:
$$1-e^{-\alpha\sqrt3}>\frac{M\sqrt3}{|R|^2}\quad(R\gg1).$$
On the other arcs $\mathbb{T}_R^1$ and $\mathbb{T}_R^3$ the proof is analogous.

Rouche's theorem implies that the number of zeros of $B_1(\lambda,l)$ and of the function $s_1(i\lambda l)-s_2(i\lambda l)$ inside the circle $\mathbb{T}_R$ is the same. And since zeros of $B_1(\lambda,l)$ are included in pairs of points symmetrical relative to the axis $i\mathbb{R}$, then hence it follows that zeros of $B_1(\lambda,l)$ lye on the axis $i\mathbb{R}$. By $\Lambda_q$, we denote the set of zeros of the function $B_1^*(\lambda,l)$ which is complexly conjugated with the set of zeros of the function $B_1(\lambda,l)$,
\begin{equation}
\Lambda_q\stackrel{\rm def}{=}\{\mu_n(q)=i\nu_n(q):B_1^*(\mu_n(q))=0;\nu_n(q)\in\mathbb{R}\}.\label{eq3.29}
\end{equation}

\begin{lemma}\label{l3.8}
Equality \eqref{eq3.25} generates the problem of jump (on the ray $il_{\zeta_2}$) in the left half-plane $i\mathbb{C}_+=\{\lambda\in\mathbb{C}:\Re\lambda<0\}$ for holomorphic in the adjacent sectors $\Omega_2$ and $\Omega_1\cap\Omega_3^-$ functions $E_2(\lambda,x)$ and $f_{1,2}(\lambda,x)$ which tend to $1$ as $\lambda\rightarrow\infty$ (inside these sectors). Function $f_{1,2}(\lambda,x)$ has no more than countable number of poles at the points of the set
\begin{equation}
\Lambda_q^1=\Lambda_q\cap\Omega_1(\in i\mathbb{R}_-)\label{eq3.30}
\end{equation}
where $\Lambda_q$ is given by \eqref{eq3.29}.
\end{lemma}

Analogously, relation (ii) \eqref{eq3.18},
\begin{equation}
\begin{array}{ccc}
\zeta_3e^{i\lambda(\zeta_3-\zeta_1)}c_3^*(\lambda\zeta_3)E_3(\lambda,x)=E_1(\lambda,x)-f_{2,3}(\lambda,x)\\
{\displaystyle\left(f_{2,3}(\lambda,x)\stackrel{\rm def}{=}\frac{\zeta_1}{\sqrt3\lambda}e^{-i\lambda\zeta_1x}\omega_{2,3}^*(\lambda,x)\right),}
\end{array}\label{eq3.31}
\end{equation}
gives the problem of jump on the ray $il_{\zeta_3}$ in the adjacent sectors $\Omega_1$ and $\Omega_3\cap\Omega_2^-$ of the half-plane $i\zeta_2\mathbb{C}_+$ for the analytical in these sectors functions $E_1(\lambda,x)$ and $f_{2,3}(\lambda,x)$ tending to $1$ as $\lambda\rightarrow\infty$ (insides the sectors). Function $f_{2,3}(\lambda,x)$ has poles at the points $\zeta_2\Lambda_q^1\in(-il_{\zeta_2})$ where $\Lambda_q^1$ is given by \eqref{eq3.30}.

Similarly, equality (iii) \eqref{eq3.18}
\begin{equation}
\begin{array}{ccc}
\zeta_3e^{i\lambda(\zeta_2-\zeta_3)x}c_3^*(\lambda\zeta_2)E_2(\lambda,x)=E_3(\lambda,x)-f_{3,1}(\lambda,x)\\
{\displaystyle\left(f_{3,1}(\lambda,x)\stackrel{\rm def}{=}\frac{\zeta_2}{\sqrt3\lambda}e^{-i\lambda\zeta_3x}\omega_{3,1}^*(\lambda,x)\right)}
\end{array}\label{eq3.32}
\end{equation}
is the problem of jump on the ray $il_{\zeta_1}$ in the half-plane $i\zeta_3\mathbb{C}_+$ for holomorphic in adjacent sectors $\Omega_3$ and $\Omega_2\cap\Omega_1^-$ functions $E_3(\lambda,x)$ and $f_{3,1}(\lambda,x)$ which tend to 1 when $\lambda\rightarrow\infty$ (inside the sectors). Moreover, $f_{3,1}(\lambda,x)$ has poles at the points of the set $\zeta_3\Lambda_q^1\in(-il_{\zeta_3})$.

Equalities \eqref{eq3.17} give three more problems of jumps in adjacent sectors of the half-plane. So, relation (i) \eqref{eq3.17}
\begin{equation}
\zeta_2e^{i\lambda(\zeta_1-\zeta_3)x}c_2^*(\lambda)E_1(\lambda,x)=E_3(\lambda,x)-g_{1,3}(\lambda,x)\quad\left(g_{1,3}(\lambda,x)\stackrel{\rm def}{=}-\frac{\zeta_2e^{-i\lambda\zeta_3x}}{\sqrt3\lambda}\omega_{1,3}^*(\lambda,x)\right)\label{eq3.33}
\end{equation}
is the problem of jump on the ray $il_{\zeta_3}$ in the right half-plane $-i\mathbb{C}_+=\{\lambda\in\mathbb{C}:\Re\lambda>0\}$ for the analytical in the adjacent sectors $\Omega_3$ and $\Omega_1\cap\Omega_2^-$ functions $E_3(\lambda,x)$ and $g_{1,3}(\lambda,x)$ that tend to 1 as $\lambda\rightarrow\infty$ (inside these sectors). Function
$g_{1,3}(\lambda,x)$ has poles at the points of $\Lambda_q^1$ \eqref{eq3.20}.

Relation (ii) \eqref{eq3.17}
\begin{equation}
\begin{array}{ccc}
\zeta_2e^{i\lambda(\zeta_3-\zeta_2)x}c_2^*(\lambda\zeta_3)E_3(\lambda,x)=E_2(\lambda,x)-g_{2,1}(\lambda,x)\\
{\displaystyle\left(g_{2,1}(\lambda,x)\stackrel{\rm def}{=}-\frac{\zeta_3}{\sqrt3\lambda}e^{-i\lambda\zeta_2x}\omega_{2,1}^*(\lambda,x)\right)}
\end{array}\label{eq3.34}
\end{equation}
gives the problem of jump on the ray $il_{\zeta_1}$ in the half-plane $(-i\zeta_2\mathbb{C}_+)$ for the holomorphic in adjacent sectors $\Omega_2$ and $\Omega_3\cap\Omega_1^-$ functions $E_2(\lambda,x)$ and $g_{2,1}(\lambda,x)$ tending to 1 as $\lambda\rightarrow\infty$ (inside the sectors). Function $g_{2,1}(\lambda,x)$ has poles at the points of the set $\zeta_2\Lambda_q^1\subset(-il_{\zeta_2})$.

Finally, equality (iii) \eqref{eq3.17}
\begin{equation}
\begin{array}{ccc}
\zeta_2e^{i\lambda(\zeta_2-\zeta_1)x}c_2^*(\lambda\zeta_2)E_2(\lambda,x)=E_1(\lambda,x)-g_{3,2}(\lambda,x)\\
{\displaystyle\left(g_{3,2}(\lambda,x)\stackrel{\rm def}{=}-\frac{\zeta_1}{\sqrt3\lambda}e^{-i\lambda\zeta_1x}\omega_{3,2}^*(\lambda,x)\right)}
\end{array}\label{eq3.35}
\end{equation}
generates the problem of jump on the ray $il_{\zeta_2}$ in the half-plane $(-i\zeta_3\mathbb{C}_+)$ for analytical in the adjacent sectors $\Omega_1$ and $\Omega_2\cap\Omega_3^-$ functions $E_1(\lambda,x)$ and $g_{3,2}(\lambda,x)$ tending to 1 for $\lambda\rightarrow\infty$ (inside the sectors). Function $g_{3,2}(\lambda,x)$ has poles at the points of the set $\zeta_3\Lambda_q^1\subset(-il_{\zeta_3})$.
\vspace{5mm}

{\bf 3.6} Consider two problems of the jump \eqref{eq3.31} and \eqref{eq3.6} on the rays $il_{\zeta_3}$ and $il_{\zeta_2}$,
\begin{equation}
\begin{array}{lll}
\zeta_3e^{i\lambda(\zeta_3-\zeta_1)}c_3^*(\lambda\zeta_3)E_3(\lambda,x)=E_1(\lambda,x)-f_{2,3}(\lambda,x)\quad(\lambda\in il_{\zeta_3});\\
\zeta_2e^{i\lambda(\zeta_2-\zeta_1)}c_2^*(\lambda\zeta_2)E_2(\lambda,x)=E_1(\lambda,x)-g_{3,2}(\lambda,x)\quad(\lambda\in il_{\zeta_2})
\end{array}\label{eq3.36}
\end{equation}
where $E_1(\lambda,x)$ is holomorphic in $\Omega_1$, and $f_{2,3}(\lambda,x)$ and $g_{3,2}(\lambda,x)$, in the sectors $\Omega_3\cap\Omega_2^-$ and $\Omega_2\cap\Omega_3^-$ correspondingly.

\begin{picture}(200,200)
\put(0,100){\vector(1,0){200}}
\put(100,0){\vector(0,1){200}}
\put(0,133){\vector(3,-1){200}}
\put(100,100){\vector(3,1){100}}
\put(100,100){\vector(-3,-1){100}}
\put(2,130){$-il_{\zeta_3}$}
\put(100,150){$g_{2,1}(\lambda,x)$}
\put(40,150){$f_{3,1}(\lambda,x)$}
\put(100,180){$il_{\zeta_1}$}
\qbezier(100,140)(135,145)(141,115)
\qbezier(140,85)(155,105)(143,116)
\qbezier(100,60)(135,70)(144,85)
\qbezier(100,58)(70,60)(57,82)
\qbezier(56,86)(50,100)(58,115)
\qbezier(100,141)(65,150)(57,116)
\put(10,55){$il_{\zeta_2}$}
\put(110,55){$E_1(\lambda,x)$}
\put(0,89){$g_{3,2}(\lambda,x)$}
\put(0,105){$\zeta_3\mu_n(q)$}
\put(180,135){$-il_{\zeta_2}$}
\put(198,105){$l_{\zeta_1}$}
\put(169,53){$il_{\zeta_3}$}
\put(169,80){$f_{2,3}(\lambda,x)$}
\put(50,112){$\circ$}
\put(150,114){$\circ$}
\put(150,105){$\zeta_2\mu_n(q)$}
\end{picture}

\hspace{20mm} Fig. 3

Using \eqref{eq3.31}, \eqref{eq3.36} and \eqref{eq3.4}, \eqref{eq3.5}, we find the boundary values of the functions $f_{2,3}(\lambda,x)$ and $g_{3,2}(\lambda,x)$ on the rays $(-il_{\zeta_2})$ and $(-il_{\zeta_3})$,
\begin{equation}
\begin{array}{lll}
{\displaystyle f_{2,3}(-i\zeta_2t,x)=\frac{i\zeta_3}{\sqrt3t}e^{-\zeta_2tx}\frac1{B_1^*(-it)}\{W^*(-it,x),e_2^*(-it,x)\}\quad(t\geq0);}\\
{\displaystyle g_{3,2}(-i\zeta_3t,x)=-\frac{i\zeta_2}{\sqrt3t}e^{-\zeta_3tx}\frac1{B_1^*(-it)}\{W^*(-it,x),e_3^*(-it,x)\}\quad(t\geq0).}
\end{array}\label{eq3.37}
\end{equation}
Boundary values on the rays $(-il_{\zeta_3})$ and $(-il_{\zeta_2})$ of holomorphic in the sectors $\Omega_2\cap\Omega_1^-$ and $\Omega_3\cap\Omega_1^-$ functions $f_{3,1}(\lambda,x)$ and $g_{2,1}(\lambda,x)$ are
\begin{equation}
\begin{array}{lll}
{\displaystyle f_{3,1}(-i\zeta_3t,x)=\frac{i\zeta_3}{\sqrt3t}e^{-\zeta_2tx}\frac1{B_1^*(-it)}\{W^*(-it,x),e_2^*(-it,x)\}\quad(t\geq0);}\\
{\displaystyle g_{2,1}(-i\zeta_2t,x)=-\frac{i\zeta_2}{\sqrt3}e^{-\zeta_3tx}\frac1{B_1^*(-it)}\{W^*(-it,x),e_3^*(-it,x)\}\quad(t\geq0).}
\end{array}\label{eq3.38}
\end{equation}
Equations \eqref{eq3.37} and \eqref{eq3.38} imply
$$\left.f_{2,3}(\lambda,x)\right|_{\lambda\in(-il_{\zeta_2})}=\left.f_{3,1}(\lambda,x)\right|_{\lambda\in(-il_{\zeta_3})};\quad\left.g_{3,2}(\lambda,x)\right|_{\lambda\in(-il_{
\zeta_3})}=\left.g_{2,1}(\lambda,x)\right|_{\lambda\in(-il_{\zeta_2})}.$$

Symmetry $\lambda\rightarrow-\overline{\lambda}$ relative to the axis $iL_{\zeta_1}$ establishes a one-to-one correspondence between the sectors $\Omega_2\cap\Omega_1^-$ and $\Omega_3\cap\Omega_1^-$, and let $f^+(\lambda)\stackrel{\rm def}{=}f(-\overline{\lambda})$. The last equalities imply that values of $f_{2,3}(\lambda,x)$ and $f_{3,1}^+(\lambda,x)$ on the ray $(-il_{\zeta_2})$ (and also the values of $g_{3,2}(\lambda,x)$ and $g_{2,1}^+(\lambda,x)$ on the ray $(-il_{\zeta_3})$) coinside and are continuous (except for singularities $\zeta_2\Lambda_q^1$ and $\zeta_3\Lambda_q^1$), therefore the functions
\begin{equation}
\begin{array}{lll}
G_2(\lambda,x)\stackrel{\rm def}{=}g_{3,2}(\lambda,x)\chi_{\Omega_2\cap\Omega_3^-}+g_{2,1}^+(\lambda,x)\chi_{\Omega_2\cap\Omega_1^-}\quad(\lambda\in\Omega_2);\\
F_3(\lambda,x)\stackrel{\rm def}{=}f_{2,3}(\lambda,x)\chi_{\Omega_3\cap\Omega_2^-}+f_{3,1}^+(\lambda,x)\chi_{\Omega_3\cap\Omega_1^-}\quad(\lambda\in\Omega_3)
\end{array}\label{eq3.39}
\end{equation}
are holomorphic in the sectors $\Omega_2$ and $\Omega_3$ (except for poles at the points of $\zeta_3\Lambda_q^1$ and $\zeta_2\Lambda_q^1$). Moreover, $G_2(\lambda,x)$ and $F_3(\lambda,x)$ tend to 1 when $\lambda\rightarrow\infty$, for $\lambda\in\Omega_2\backslash(-il_{\zeta_3})$ and $\lambda\in\Omega_3\backslash(-il_{\zeta_2})$. Boundary values of $f_{3,1}(\lambda,x)$ and $g_{2,1}(\lambda,x)$ on the ray $(il_{\zeta_1})$ are
$$f_{3,1}(it,x)(=f_{3,1}^+(it,x))=-\frac{i\zeta_2}{\sqrt3t}e^{\zeta_3tx}\frac{\{W^*(it,x),e_1^*(it,x)\}}{B_3^*(it)};$$
$$g_{2,1}(it,x)(=g_{2,1}^+(it,x))=\frac{i\zeta_3}{\sqrt3t}e^{\zeta_2tx}\frac{\{W^*(it,x),e_1^*(it,x)\}}{B_2^*(it)},$$
therefore
\begin{equation}
G_2(it,x)=d(it,x)F_3(it,x)\quad(d(\lambda,x)\stackrel{\rm def}{=}-\zeta_2e^{-\sqrt3\lambda x}c_2^*(\lambda\zeta_3)).\label{eq3.40}
\end{equation}

As a result, we arrive at the Riemann boundary value problem \cite{20,21} on the contour $\Gamma\stackrel{\rm def}{=}\bigcup\limits_k(il_{\zeta_k})$ formed by the rays $il_{\zeta_k}$ ($1\leq k\leq3$).

\begin{picture}(200,200)
\put(0,100){\vector(1,0){200}}
\put(100,0){\vector(0,1){200}}
\put(95,200){\vector(0,-1){50}}
\put(95,150){\line(0,-1){45}}
\put(95,105){\vector(-3,-1){50}}
\put(50,90){\line(-3,-1){50}}
\put(0,60){\vector(3,1){50}}
\put(50,77){\line(3,1){50}}
\put(100,94){\vector(3,-1){50}}
\put(150,77){\line(3,-1){50}}
\put(200,73){\vector(-3,1){50}}
\put(150,89){\line(-3,1){45}}
\put(105,104){\vector(0,1){50}}
\put(105,154){\line(0,1){45}}
\put(100,100){\vector(-3,-1){100}}
\put(100,100){\vector(3,-1){100}}
\put(180,55){$il_{\zeta_3}$}
\put(120,150){$F_3(\lambda,x)$}
\put(10,55){$il_{\zeta_2}$}
\put(40,150){$G_2(\lambda,x)$}
\put(105,180){$il_{\zeta_1}$}
\put(105,35){$E_1(\lambda,x)$}
\put(180,105){$l_{\zeta_1}$}
\end{picture}

\hspace{20mm} Fig. 4

On the rays $il_{\zeta_2}$ and $il_{\zeta_3}$, the boundary value problem coincides with the problem of jump \eqref{eq3.36} where $f_{2,3}(\lambda,x)=F_3(\lambda,x)$ for $\lambda\in(il_{\zeta_3})$ and $g_{3,2}(\lambda,x)=G_2(\lambda,x)$ for $\lambda\in(il_{\zeta_2})$. On the ray $(il_{\zeta_1})$, boundary condition is given by \eqref{eq3.40}. Solution to such a boundary value problem is realized in a standard manner \cite{20,21}. By $\chi(\lambda,x)$, we define the {\bf canonical solution} \cite{20,21}
\begin{equation}
\chi(\lambda,x)\stackrel{\rm def}{=}\exp\left\{\frac1{2\pi i}\int\limits_\Gamma\frac{\ln D(\tau,x)}{\tau-\lambda}\right\}=\exp\left\{\frac1{2\pi i}\int\limits_0^\infty\frac{\ln d(i\tau,x)}{\tau+i\lambda}\right\}\label{eq3.41}
\end{equation}
and define the piecewise holomorphic function in $\mathbb{C}$
\begin{equation}
\mathcal{F}(\lambda,x)\stackrel{\rm def}{=}\left\{
\begin{array}{ccc}
E_1(\lambda,x)\chi^{-1}(\lambda,x)\quad(\lambda\in\Omega_1);\\
G_2(\lambda,x)\chi^{-1}(\lambda,x)\quad(\lambda\in\Omega_2);\\
F_3(\lambda,x)\chi^{-1}(\lambda,x)\quad(\lambda\in\Omega_3);
\end{array}\right.\label{eq3.42}
\end{equation}
which is analytical in the sectors $\{\Omega_k\}_1^3$ (except for poles at $\zeta_2\Lambda_q^1$ and at $\zeta_3\Lambda_q^1$) and tends to 1 when $\lambda\rightarrow\infty$. Function $\mathcal{F}(\lambda,x)$ is holomorphic on the ray $il_{\zeta_1}$, and the boundary value problem is reduced to the problem of jump on the broken line $\widetilde{\Gamma}\stackrel{\rm def}{=}(il_{\zeta_2})\cup(il_{\zeta_3})$, besides, $(il_{\zeta_2})$ is an incoming ray and $(il_{\zeta_3})$ is the outgoing ray (see Fig. 4). This problem of a jump follows from \eqref{eq3.36} after division by $\chi(\lambda,x)$ \eqref{eq3.41}. Hence it follows that
\begin{equation}
\mathcal{F}(\lambda,x)=1+b(\lambda,x)+\frac1{2\pi i}\int\limits_{\widetilde{\Gamma}}\frac{\varphi(\tau,x)}{\tau-\lambda}d\tau\label{eq3.43}
\end{equation}
where $\varphi(\tau,x)$ is the jump function on each of the rays of the contour $\widetilde{\Gamma}$ and $b(\lambda,x)$ is given by
\begin{equation}
b(\lambda,x)\stackrel{\rm def}{=}\sum\limits_n\frac{r_n(x)}{\lambda-\zeta_2\mu_n(q)}+\sum\limits_n\frac{p_n(x)}{\lambda-\zeta_3\mu_n(q)},\label{eq3.44}
\end{equation}
besides, $\mu_n(q)\in\Lambda_q^1$ \eqref{eq3.30} and $r_n(x)$ ($p_n(x)$) is a residue of the function $F_3(\lambda,x)\chi^{-1}(\lambda,x)$ ($G_2(\lambda,x)\chi^{-1}(\lambda,x)$) at the point $\zeta_2\mu_n(q)$ ($\zeta_3\mu_n(q)$). The integral in \eqref{eq3.43} along the ray $(il_{\zeta_3})$ ($\tau=i\zeta_3t$, $t\geq0$) of the contour $\widetilde{\Gamma}$ outgoing from the origin equals
$$\frac1{2\pi i}\int\limits_{il_{\zeta_3}}\frac{\varphi(\tau,x)}{\tau-\lambda}d\tau=\frac1{2\pi i}\int\limits_0^\infty\frac{i\zeta_3dt}{i\zeta_3t-\lambda}\zeta_3e^{-(\zeta_2-\zeta_3)tx}c_3^*(i\zeta_2t)\chi^{-1}(i\zeta_3t,x)E_3(i\zeta_3t,x)$$
$$=\frac{\zeta_3}{2\pi i}\int\limits_0^\infty\frac{dt}{t+i\zeta_2\lambda}e^{-i\sqrt3tx}c_3^*(i\zeta_2t)\chi^{-1}(i\zeta_3t,x)E_2(it,x).$$
Analogously, along the incoming ray $(il_{\zeta_2})$,
$$\frac1{2\pi i}\int\limits_{il_{\zeta_2}}\frac{\varphi(\tau,x)}{\tau-\lambda}d\tau=\frac{\zeta_2}{2\pi i}\int\limits_{-\infty}^0\frac{dt}{t+i\zeta_3\lambda}e^{i\sqrt3tx}c_2^*(i\zeta_2t)\chi^{-1}(i\zeta_2t)E_3(it,x).$$
Upon substituting these expressions into \eqref{eq3.43}, we obtain
\begin{equation}
\begin{array}{ccc}
{\displaystyle\mathcal{F}(\lambda,x)=1+b(\lambda,x)+\frac{\zeta_3}{2\pi i}\int\limits_0^\infty\frac{d\tau}{\tau+i\zeta_2\lambda}D_3(i\zeta_2\tau,x)E_2(i\tau,x)}\\
{\displaystyle-\frac{\zeta_2}{2\pi i}\int\limits_0^\infty\frac{d\tau}{\tau+i\zeta_3\lambda}D_2(i\zeta_3\tau,x)E_3(i\tau,x)}
\end{array}\label{eq3.45}
\end{equation}
where $D_3(\lambda,x)$ and $D_2(\lambda,x)$ depend only on $c_2(\lambda)$, $c_3(\lambda)$,
\begin{equation}
D_3(\lambda,x)\stackrel{\rm def}{=}c_2^*(\lambda)\chi^{-1}(\lambda\zeta_2,x)e^{-\sqrt3\zeta_3\lambda x};\quad D_2(\lambda,x)\stackrel{\rm def}{=}c_3^*(\lambda)\chi^{-1}(\lambda\zeta_3,x)e^{\sqrt3\zeta_2\lambda x},\label{eq3.46}
\end{equation}
besides, $\chi(\lambda,x)$ is given by \eqref{eq3.41} and $d(\lambda,x)$, by \eqref{eq3.40} correspondingly. Assuming in \eqref{eq3.45} that $\lambda\in\Omega_1$ and using \eqref{eq3.42}, we have
\begin{equation}
\begin{array}{ccc}
{\displaystyle E_1(\lambda,x)\chi^{-1}(\lambda,x)=1+b(\lambda,x)+\frac{\zeta_3}{2\pi i}\int\limits_0^\infty\frac{d\tau}{\tau+i\zeta_2\lambda}D_3(i\zeta_2\tau,x)E_2(i\tau,x)}\\
{\displaystyle-\frac{\zeta_2}{2\pi i}\int\limits_0^\infty\frac{d\tau}{\tau+i\zeta_3\lambda}D_2(i\zeta_3\tau,x)E_3(i\tau,x).}
\end{array}\label{eq3.47}
\end{equation}
So, the holomorphic in $\Omega_1$ function $E_1(\lambda,x)\chi^{-1}(\lambda,x)$ is expressed via the functions $E_2(i\tau,x)$, $E_3(i\tau,x)$, and $b(\lambda,x)$. Calculating the boundary values in equality \eqref{eq3.47} when $\lambda\rightarrow i\zeta_2t$ and $\lambda\rightarrow i\zeta_3t$ ($\lambda\in\Omega_1$), we obtain the system of integral equations relative to $E_2(it,x)$, $E_3(it,x)$,
\begin{equation}
\left\{
\begin{array}{lclc}
{\displaystyle E_2(it,x)\chi^{-1}(it\zeta_2,x)=1+b(it\zeta_2,x)+\frac{\zeta_3}{2\pi i}\int\limits_0^\infty\frac{d\tau}{\tau-\zeta_3t}D_3(i\zeta_2\tau,x)E_2(i\tau,x)}\\
{\displaystyle-\frac{\zeta_2}{2\pi i}\int\limits_0^\infty\hspace{-4.4mm}/\frac{d\tau}{\tau-t}D_2(i\zeta_3\tau,x)E_3(i\tau,x)-\frac{\zeta_2}2D_2(i\zeta_3t,x)E_3(it,x);}\\
{\displaystyle E_3(it,x)\chi^{-1}(it\zeta_3,x)=1+b(it\zeta_3,x)+\frac{\zeta_3}{2\pi i}\int\limits_0^\infty\hspace{-4.4mm}/\frac{d\tau}{\tau-t}D_3(i\zeta_2\tau,x)E_2(i\tau,x)}\\
{\displaystyle-\frac{\zeta_2}{2\pi i}\int\limits_0^\infty\frac{d\tau}{\tau-\zeta_2t}D_2(i\zeta_3\tau,x)E_3(i\tau,x)+\frac{\zeta_3}2D_3(i\zeta_2t)E_2(it,x).}
\end{array}\right.\label{eq3.48}
\end{equation}
Also, we need the relations to find $r_n(x)$ and $p_n(x)$ from \eqref{eq3.44}. Equations \eqref{eq3.31} and \eqref{eq3.35} imply
$$f_{2,3}(\lambda\zeta_2,x)=\frac{\zeta_3}{\sqrt3\lambda}e^{-i\lambda\zeta_2x}\frac{B_3^*(\lambda)}{B_2^*(\lambda\zeta_2)}\omega_{3,2}^*(\lambda,x)=-\zeta_3e^{i\lambda(\zeta_1-
\zeta_2)x}\frac{B_3^*(\lambda)}{B_1^*(\lambda)}g_{3,2}(\lambda,x),$$
and thus
$$\frac{f_{2,3}(\lambda\zeta_2,x)}{B_3^*(\lambda)}=-\zeta_3e^{-\sqrt3\zeta_3\lambda x}\frac{g_{3,2}(\lambda,x)}{B_1^*(\lambda,x)}.$$
Since $f_{2,3}(\lambda,x)\chi^{-1}(\lambda,x)=\mathcal{F}(\lambda,x)$ for $\lambda\in\Omega_3\cap\Omega_2^-$ ($g_{3,2}(\lambda,x)\chi^{-1}(\lambda,x)=\mathcal{F}(\lambda,x)$ for $\lambda\in\Omega_2\cap\Omega_3^-$), then the equality
$$\frac{f_{2,3}(\lambda\zeta_2,x)}{B_3^*(\lambda)}\chi^{-1}(\lambda\zeta_2,x)=-\zeta_3e^{-\sqrt3\zeta_3\lambda x}\chi^{-1}(\lambda\zeta_2,x)\chi(\lambda,x)\frac{g_{3,2}(\lambda,x)}{B_1^*(\lambda,x)}\chi^{-1}(\lambda,x),$$
upon calculation of the residue at the point $\lambda=\zeta_3\mu_n(q)$, we obtain
$$\frac1{\dot{B}_1(\mu_n(q))}\left\{1+\sum\limits_n\frac{\zeta_3r_n(x)}{\zeta_2\mu_m(q)-\mu_n(q)}+\sum\limits_{n\not=m}\frac{\zeta_2p_m(x)}{\mu_m(q)-\mu_n(q)}+\frac{
\zeta_3}{2\pi i}\int\limits_0^\infty\frac{d\tau}{\tau+i\mu_m(q)}\right.$$
\begin{equation}
\left.\times D_3(i\zeta_2\tau,x)E_2(i\tau,x)-\frac{\zeta_2}{2\pi i}\int\limits_0^\infty\frac{d\tau}{\tau+i\zeta_2\mu_m(q)}D_2(i\zeta_3\tau,x)E_3(i\tau,x)\right\}\label{eq3.49}
\end{equation}
$$=-\zeta_3e^{-\sqrt3\zeta_2\mu_m(q)x}\chi^{-1}(\mu_m(q),x)\chi(\zeta_3\mu_m(q),x)\frac{p_n(x)}{B_2^*(\mu_m(q))}\quad(m\in\mathbb{Z}_+)$$
($\dot{B}_1(\lambda)$ is derivative with respect to $\lambda$). Analogously, taking into account
$$\frac{g_{3,2}(\lambda\zeta_3,x)}{B_2^*(\lambda)}=-\zeta_2e^{\sqrt3\zeta_2x}\frac{f_{2,3}(\lambda,x)}{B_1^*(\lambda)},$$
calculate the residue at the point $\lambda=\zeta_2\mu_m(q)$, then
$$\frac1{\dot{B}_1(\mu_m(q))}\left\{1+\sum\limits_{n\not=m}\frac{\zeta_3r_n(x)}{\mu_m(q)-\mu_n(q)}+\sum\limits_n\frac{\zeta_2p_n(x)}{\zeta_3\mu_m(q)-\mu_n(q)}+\frac{\zeta_3}{2\pi i}\int\limits_0^\infty\frac{d\tau}{\tau+i\zeta_3\mu_m(q)}\right.$$
\begin{equation}
\left.\times D_3(i\zeta_2\tau,x)E_2(i\tau,x)-\frac{\zeta_2}{2\pi i}\int\limits_0^\infty\frac{d\tau}{\tau+i\mu_m(q)}D_2(i\zeta_3\tau,x)E_3(i\tau,x)\right\}\label{eq3.50}
\end{equation}
$$=-\zeta_2e^{\sqrt3\zeta_3\mu_m(q)x}\chi^{-1}(\mu_m(q),x)\chi(\zeta_2\mu_m(q),x)\frac{r_m(x)}{B_3^*(\mu_m(q))}\quad(m\in\mathbb{Z}_+).$$
Notice that $\mu_n(q)=i\nu_n(q)$ \eqref{eq3.29} and $\mu_n(q)\in i\mathbb{R}_-$ \eqref{eq3.30}, i. e., $\nu_n(q)<0$, and consequently, integrals in \eqref{eq3.49}, \eqref{eq3.50} don't have singularities.

{\bf Conclusion.} {\it Relations \eqref{eq3.48} -- \eqref{eq3.50} ($b(\lambda,x)$ is given by \eqref{eq3.44}) form a {\bf closed system of linear singular equations} relative to $E_2(it,x)$, $E_3(it,x)$ and $\{r_n(x)\}$, $\{p_n(x)\}$.}

\begin{remark}\label{r3.3}
Function $B_1(\lambda)$ is the independent parameter of the system \eqref{eq3.48} -- \eqref{eq3.50}. From $B_1(\lambda)$, we find $B_2(\lambda)=B_1(\lambda\zeta_2)$ and $B_3(\lambda)=B_1(\lambda\zeta_3)$ (see \eqref{eq3.4}) and define $c_2(\lambda)$, $c_3(\lambda)$ \eqref{eq3.7}. The numbers $\{\mu_n(q)\}$, zeros of $B_1^*(\lambda)$, and constants
\begin{equation}
a_m\stackrel{\rm def}{=}\dot{B}_1(\mu_m(q))/B_3(\mu_m(q));\quad b_m\stackrel{\rm def}{=}\dot{B}_1(\mu_m(q))/B_2(\mu_m(q))\quad(m\in\mathbb{R}_+)\label{eq3.51}
\end{equation}
are calculated from the function $B_1(\lambda)$.
\end{remark}

\section{Inverse problem}\label{s4}

{\bf 4.1} Method of solution of the inverse problem follows from considerations of Section 3. Knowing the functions $s_1(\lambda,l)$ and $s_2(\lambda,l)$ (or $\widehat{s}_p(\lambda,l)=(i\lambda)^ps_p(\lambda,l)$; $p=1$, 2), we define the functions $\{B_k(\lambda)\}_1^3$ \eqref{eq3.4}, and then $\{c_k(\lambda)\}_2^3$ \eqref{eq3.5} also, and, finally, we define $\Lambda_q^1$ \eqref{eq3.30}, the set of zeros of $B_1^*(\lambda)$ situated on $i\mathbb{R}_-$. From the set $\{c_2(\lambda),c_3(\lambda),\Lambda_q^1\}$, we construct the system of linear equations \eqref{eq3.48} -- \eqref{eq3.50}. Upon solving this system, we find $E_2(it,x)$, $E_3(it,x)$ and $\{r_n(x)\}$, $\{p_n(x)\}$, and define the function $E_1(\lambda,x)$ \eqref{eq3.47}, and then, using (a) \eqref{eq3.22}, we calculate
\begin{equation}
\lim\limits_{\lambda\rightarrow\infty}3i\lambda^2\{E_1(\lambda,x)-1\}=\int\limits_0^xq(t)dt\quad(\lambda\in\Omega_1),\label{eq4.1}
\end{equation}
hence, upon differentiating, we obtain potential $q(x)$.

The scheme stated above is based upon the solution of two problems: (i) finding the functions $s_1(\lambda,l)$, $s_2(\lambda,l)$ from the {\bf spectral data}; (ii) solvability of the system of linear equations \eqref{eq3.48} -- \eqref{eq3.50}.
\vspace{5mm}

{\bf 4.2} From the set of numbers
\begin{equation}
{\mathfrak A}_\theta\stackrel{\rm def}{=}(a;\theta;\{\lambda_n(q,\theta)\})\label{eq4.2}
\end{equation}
where $a\in\mathbb{C}$, ($a\not=0$); $\theta\in\mathbb{T}$ ($\theta+\theta_0\not=0$, $\theta_0=\overline{a}/a$) \eqref{eq2.37}, we construct the function $\Delta_\theta(q,\lambda)$ \eqref{eq2.37}, due to Lemma \ref{l2.6}, and using \eqref{eq2.34} we obtain
\begin{equation}
\frac1{\theta+\theta_0}(\theta s_2(\lambda,l)+s_2^*(\lambda,l))=\Pi(\theta)\quad(\Pi(\theta)\stackrel{\rm def}{=}a\Pi\left(1-\frac{\lambda^3}{\lambda_n^3(q,\theta)}\right)).\label{eq4.3}
\end{equation}
Analogously, relation
\begin{equation}
\frac1{\widehat{\theta}+\theta_0}(\widehat{\theta}s_2(\lambda,l)+s_2^*(\lambda,l))=\Pi(\widehat{\theta})\label{eq4.4}
\end{equation}
corresponds to another set $\mathfrak{A}_{\widehat{\theta}}$ ($\widehat{\theta}\not=\theta$, $\widehat{\theta}+\theta_0\not=0$). Solution to the system \eqref{eq4.3}, \eqref{eq4.4} is given by
\begin{equation}
s_2(\lambda,l)=\frac1{\theta-\widehat{\theta}}\{\Pi(\theta)(\theta+\theta_0)-\Pi(\widehat{\theta})(\widehat{\theta}+\theta_0)\}.\label{eq4.5}
\end{equation}
So, the function $s_2(\lambda,l)$ \eqref{eq4.5} (and thus $s_2^*(\lambda,l)$ also) is calculated unambiguously from the two spectral sets $\mathfrak{A}_\theta$ and $\mathfrak{A}_{\widehat{\theta}}$ \eqref{eq4.2} ($\theta\not=\widehat{\theta}$, $\theta+\theta_0\not=0$).

By $L_q(\theta,h)$, we denote the self-adjoint operator in $L^2(0,l)$,
\begin{equation}
(L_q(\theta,h)y)(x)\stackrel{\rm def}{=}y'''(x)+q(x)y(x)\label{eq4.6}
\end{equation}
($q(x)$ is a real function from $L^2(0,l)$), domain of which is given by
\begin{equation}
\mathfrak{D}(L_q(\theta,h))\stackrel{\rm def}{=}\{y\in W_2^3(0,l):y(0)=ihy'(0);y'(l)=\theta y'(0);y(l)=0\}\label{eq4.7}
\end{equation}
($h\in\mathbb{R}$; $\theta\in\mathbb{T}$). The function
$$Y(\lambda,x)\stackrel{\rm def}{=}y_0s_0(\lambda,x)+y_1s_1(\lambda,x)+y_2s_2(\lambda,x)$$
($y_k\in\mathbb{C}$; $s_k(\lambda,x)$ is given by \eqref{eq2.25}) is the solution to equation \eqref{eq2.1} and boundary conditions \eqref{eq4.7} for $Y(\lambda,x)$ give the equation system
$$\left\{
\begin{array}{lll}
y_0-ihy_2=0;\\
y_0s'_0(\lambda,l)+y_1(s'_1(\lambda,l)-\theta)+y_2s'_2(\lambda,l)=0;\\
y_0s_0(\lambda,l)+y_1s_1(\lambda,l)+y_2s_2(\lambda,l)=0.
\end{array}\right.$$
This system has a non-trivial solution if its determinant $\Delta_{\theta,h}(q,\lambda)$ vanishes, where
$$\Delta_{\theta,h}(q,\lambda)\stackrel{\rm def}{=}\det\left[
\begin{array}{ccc}
1&0&-ih\\
s'_0(\lambda,l)&s'_1(\lambda,l)-\theta&s'_2(\lambda,l)\\
s_0(\lambda,l)&s_1(\lambda,l)&s_2(\lambda,l)
\end{array}\right]$$
$$=(s'_1(\lambda,l)-\theta)s_2(\lambda,l)-s'_2(\lambda,l)s_1(\lambda,l)-ih[s'_0(\lambda,l)s_0(\lambda,l)-(s'_1(\lambda,l)-\theta)s_0(\lambda,l)].$$
Taking into account \eqref{eq2.32}, we arrive at the characteristic function of the operator $L_q(\theta,h)$
\begin{equation}
\Delta_{\theta,h}(q,\lambda)=-\theta s_2(\lambda,l)-s_2^*(\lambda,l)-ih(\theta s_2(\lambda,l)-s_0^*(\lambda,l)),\label{eq4.8}
\end{equation}
and, due to \eqref{eq2.34},
\begin{equation}
\Delta_{\theta,h}(q,\lambda)=\Delta_\theta(q,\lambda)-ih(\theta s_0(\lambda,l)-s_0^*(\lambda,l)).\label{eq4.9}
\end{equation}
For $\Delta_{\theta,h}(q,\lambda)$ \eqref{eq4.8}, analogously to \eqref{eq2.37}, the following multiplicative decomposition holds:
\begin{equation}
\Delta_{\theta,h}(q,\lambda)=-\{a(\theta+\theta_0)+ihb(\theta-\theta_1)\}\prod\limits_n\left(1-\frac{\lambda^3}{\lambda_n^3(q,\theta,h)}\right)\label{eq4.10}
\end{equation}
where $a=s_2(\lambda,l)$; $b=s_0(\lambda,l)$; $\theta_0=\overline{a}/a$; $\theta_1=\overline{b}/b$; $\lambda_n(q,\theta,h)$ are real zeros of the function $\Delta_{\theta,h}(q,\lambda)$ \eqref{eq4.8}. Hereinafter, we assume that
\begin{equation}
a\not=0;\quad\theta+\theta_0\not=0;\quad a(\theta+\theta_0)+ihb(\theta-\theta_1)\not=0.\label{eq4.11}
\end{equation}
Conditions \eqref{eq4.11} mean that $\lambda=0$ is not a point of the spectrum of the operators $L_q(\theta)$ and $L_q(\theta,h)$. Equation \eqref{eq4.10} implies that the set
\begin{equation}
\mathfrak{A}_{\theta,h}\stackrel{\rm def}{=}(a;b;\theta,h,\{\lambda_n(q,\theta,h)\}),\label{eq4.12}
\end{equation}
for which \eqref{eq4.11} takes place, unambiguously defines the characteristic function $\Delta_{\theta,h}(q,\lambda)$ \eqref{eq4.8}.

From the sets $\mathfrak{A}_\theta$ \eqref{eq4.2}, $\mathfrak{A}_{\theta,h}$ \eqref{eq4.12}, we construct the functions $\Delta_\theta(q,\lambda)$ \eqref{eq2.37} and $\Delta_{\theta,h}(q,\lambda)$ \eqref{eq4.10}, and then, using \eqref{eq4.9}, we find
\begin{equation}
ih(\theta s_0(\lambda,l)-s_0^*(\lambda,l))=\Delta_\theta(q,\lambda)-\Delta_{\theta,h}(q,\lambda).\label{eq4.13}
\end{equation}
Analogously, from the two sets $\mathfrak{A}_\theta$ \eqref{eq4.2}, $\mathfrak{A}_{\theta,h}$ \eqref{eq4.12}, for which \eqref{eq4.11} holds, we obtain the second equation
\begin{equation}
ih(\widehat{\theta}s_0(\lambda,l)-s_0^*(\lambda,l))=\Delta_{\widehat{\theta}}(q,\lambda)-\Delta_{\widehat{\theta},h}(q,\lambda).\label{eq4.14}
\end{equation}
For $h\not=0$, from \eqref{eq4.13}, \eqref{eq4.14} we find
\begin{equation}
s_0(\lambda,l)=\frac1{ih(\theta-\widehat{\theta})}\{\Delta_\theta(q,\lambda)-\Delta_{\widehat{\theta}}(q,\lambda)-\Delta_{\theta,h}(q,\lambda)+\Delta_{\widehat{\theta},h}(q,
\lambda)\}.\label{eq4.15}
\end{equation}

\begin{lemma}\label{l4.1}
From the four spectral sets $\mathfrak{A}_\theta$, $\mathfrak{A}_{\widehat{\theta}}$ \eqref{eq4.2} ($\theta\not=\widehat{\theta}$) and $\mathfrak{A}_{\theta,h}$, $\mathfrak{A}_{\widehat{\theta},h}$ \eqref{eq4.12} ($h\not=0$) satisfying conditions \eqref{eq4.11}, the functions $s_2(\lambda,l)$ \eqref{eq4.5} and $s_0(\lambda,l)$ \eqref{eq4.15} are defined unambiguously.
\end{lemma}

Equation \eqref{eq2.32} implies
$$s_2^*(\lambda,l)=s_1(\lambda,l)s'_2(\lambda,l)-s_2(\lambda,l)s'_1(\lambda,l);\quad s_0^*(\lambda,l)=s_0(\lambda,l)s'_1(\lambda,l)-s_1(\lambda,l)s'_0(\lambda,l).$$
Upon multiplying the first equality by $s_0(\lambda,l)$ and the second, by $s_2(\lambda,l)$, and adding, we obtain
$$s_2^*(\lambda,l)s_0(\lambda,l)+s_0^*(\lambda,l)s_2(\lambda,l)=s_1(\lambda,l)[s_0(\lambda,l)s'_2(\lambda,l)-s_2(\lambda,l)s'_0(\lambda,l)]$$
and, due to \eqref{eq2.32},
\begin{equation}
s_2^*(\lambda,l)s_0(\lambda,l)+s_0^*(\lambda,l)s_2(\lambda,l)=s_1^*(\lambda,l)s_1(\lambda,l).\label{eq4.16}
\end{equation}

Show that from this equality the function $s_1(\lambda,l)$ is unambiguously found. The function $s_1(\lambda,l)$ does not have zeros in the sectors $\{\Omega_k\}$ when $|\lambda|\gg1$, in view of \eqref{eq3.23}. Zeros of the function $s_1(\lambda,l)$ lie on the pencil of rays ${\displaystyle\Gamma=\bigcup\limits_k(il_{\zeta_k})}$, except, probably, for a finite number. Analogously, zeros of $s_1^*(\lambda,l)$ are situated on the pencil of rays ${\displaystyle\Gamma^*=\bigcup\limits_k(-il_{\zeta_k})}$ (except for a finite number of points). So, sets of zeros of the functions $s_1(\lambda,x)$ and $s_1^*(\lambda,x)$ don't intersect, except for a finite number $\lambda\in\mathbb{R}$. Every zero $\mu$ of the function $s_1(\lambda,l)$ contributes the series $\{\mu\zeta_2^l\}$ ($l=0$, 1, 2) due to $s_1(\lambda\zeta_2,l)=s_1(\lambda,l)$. Using the Hadamard theorem (Theorem \ref{t1.1}), we have that
$$s_1(\lambda,l)=ae^{b\lambda}\prod\limits_n\left(1-\frac{\lambda^3}{\mu_n^3}\right)\quad(a,b\in\mathbb{C}),$$
besides, $b=0$ ($s_2(\lambda\zeta_2,l)=s_1(\lambda,l)$, $\forall\lambda$, and the number is calculated according to \eqref{eq2.26}. So,
\begin{equation}
s_1^*(\lambda,l)s_1(\lambda,l)=|a|^2\prod\limits_n\left(1-\frac{\lambda^3}{\mu_n^3}\right)\left(1-\frac{\lambda^3}{\overline{\mu}_n^3}\right),\label{eq4.17}
\end{equation}
besides, $\mu_n\in\Gamma$ and $\overline{\mu}_n\in\Gamma^*$ for $n\gg1$. Thus, the function $s_1(\lambda,l)$ is determined from \eqref{eq4.17} up to a constant from $\mathbb{T}$ which is calculated due to \eqref{eq2.26}.

\begin{lemma}\label{l4.2}
Function $s_1(\lambda,l)$ is unambiguously found from the functions $s_0(\lambda,l)$ and $s_2(\lambda,l)$ via equality \eqref{eq4.16}.
\end{lemma}
\vspace{5mm}

{\bf 4.3} Proceed to the solvability of the system \eqref{eq3.48} -- \eqref{eq3.50}. Consider the vector-valued functions
\begin{equation}
\begin{array}{cccc}
\vec{r}(x)\stackrel{\rm def}{=}\col(r_1(x),r_2(x),...);\quad\vec{p}(x)\stackrel{\rm def}{=}\col(p_1(x),p_2(x),...);\\
{\displaystyle\vec{e}\stackrel{\rm def}{=}\col(1,1,...);\quad
\vec{h}(\tau)\stackrel{\rm def}{=}\frac1{2\pi i}\col\left(\frac1{\tau+i\mu_1(q)},\frac1{\tau+i\mu_2(q)},...\right);}\\
{\displaystyle\vec{g}(\tau,\zeta_k)\stackrel{\rm def}{=}\frac1{2\pi i}\col\left(\frac1{\tau+i\zeta_k\mu_1(q)},\frac1{\tau+i\zeta_k\mu_2(q)},...\right)}
\end{array}\label{eq4.18}
\end{equation}
($k=2$, 3). In terms of these notations, the systems \eqref{eq3.49}, \eqref{eq3.50} become
\begin{equation}
\left\{
\begin{array}{cccc}
{\displaystyle A\vec{r}(x)+B(x)\vec{p}(x)=-\zeta_3\int\limits_0^\infty d\tau\vec{h}(\tau)D_3(i\zeta_2t,x)E_2(it,x)-\vec{e}}\\
{\displaystyle+\zeta_2\int\limits_0^\infty d\tau\vec{g}(\tau,\zeta_2)D_2(i\zeta_3\tau,x)E_3(i\tau,x);}\\
{\displaystyle C(x)\vec{r}(x)+D\vec{p}(x)=-\zeta_3\int\limits_0^\infty d\tau g(\tau,\zeta_3)D_3(i\zeta_2\tau,x)E_2(i\tau,x)-\vec{e}}\\
{\displaystyle+\zeta_2\int\limits_0^\infty d\tau\vec{h}(\tau)D_2(i\zeta_3\tau,x)E_3(i\tau,x)}
\end{array}\right.\label{eq4.19}
\end{equation}
where $A$, $B(x)$, $C(x)$, $D$ are matrices,
\begin{equation}
\begin{array}{ccc}
{\displaystyle A\stackrel{\rm def}{=}\left[\frac{\zeta_3}{\zeta_2\mu_m(q)-\mu_n(q)}\right];\quad D\stackrel{\rm def}{=}\left[\frac{\zeta_2}{\zeta_3\mu_m(q)-\mu_n(q)}\right];\quad B(x)\stackrel{\rm def}{=}\left[b_n(x)\delta_{m,n}\right.}\\
{\displaystyle\left.+(1-\delta_{m,n})\frac{\zeta_2}{\mu_m(q)-\mu_n(q)}\right];\quad C(x)\stackrel{\rm def}{=}\left[c_n(x)\delta_{m,n}+(1-\delta_{m,n})\frac{\zeta_3}{\mu_m(q)-\mu_n(q)}\right],}
\end{array}\label{eq4.20}
\end{equation}
besides,
\begin{equation}
b_n(x)\stackrel{\rm def}{=}\zeta_3e^{-\sqrt3\zeta_2\mu_n(q)x}\frac{\chi(\zeta_3\mu_n(q),x)}{\chi(\mu_n(q),x)}b_n^*;\quad c_n(x)\stackrel{\rm def}{=}\zeta_2e^{\sqrt3\zeta_3\mu_n(q)x}\frac{\chi(\zeta_2\mu_n(q),x)}{\chi(\mu_n(q),x)}c_n^*\label{eq4.21}
\end{equation}
($b_n$, $c_n$ are given by the formulas \eqref{eq3.51}). Solution to this system $\vec{r}(x)$ and $\vec{p}(x)$ is linearly expressed via $E_2(i\tau,x)$, $E_3(i\tau,x)$, up to an additive summand, therefore function $b(\lambda,x)$ \eqref{eq3.44} is given by
\begin{equation}
b(\lambda,x)=a(\lambda,x)+\int\limits_0^\infty d\tau b_2(\tau,\lambda,x)E_2(i\tau,x)+\int\limits_0^\infty d\tau b_3(\tau,\lambda,x)E_3(i\tau,x)\label{eq4.22}
\end{equation}
($\lambda\not=\zeta_k\mu_n(q)$, $k=2$, 3).

Using \eqref{eq4.22}, rewrite system \eqref{eq3.48} in the matrix form, and let
\begin{equation}
\vec{\varepsilon}(t,x)=(E_2(it,x),E_3(it,x)),\label{eq4.23}
\end{equation}
then
\begin{equation}
\begin{array}{ccc}
{\displaystyle\vec{\varepsilon}(t,x)Q(t,x)+\frac1{2\pi i}\int\limits_0^\infty\hspace{-4.4mm}/\frac{d\tau}{\tau-t}\vec{\varepsilon}(\tau,x)D(\tau,x)+\frac1{2\pi i}\int\limits_0^\infty d\tau\vec{\varepsilon}(\tau,x)C(\tau,t)\widehat{D}(\tau,x)}\\
{\displaystyle=\vec{A}(t,x)+\int\limits_0^\infty d\tau\vec{\varepsilon}(t,x)\vec{B}(\tau,t,x)}
\end{array}\label{eq4.24}
\end{equation}
where
$$Q(t,x)\stackrel{\rm def}{=}\left[
\begin{array}{ccc}
\chi^{-1}(i\zeta_2t,x)&{\displaystyle-\frac{\zeta_3}2D_3(i\zeta_2t,x)}\\
{\displaystyle\frac{\zeta_2}2D_2(i\zeta_3t,x)}&\chi^{-1}(i\zeta_3t,x)
\end{array}\right];\quad C(\tau,t)\stackrel{\rm def}{=}\left[
\begin{array}{ccc}
{\displaystyle\frac{-\zeta_3}{\tau-\zeta_3t}}&0\\
0&{\displaystyle\frac{\zeta_2}{\tau-\zeta_2t}}
\end{array}\right];$$
\begin{equation}
D(\tau,x)\stackrel{\rm def}{=}\left[
\begin{array}{ccc}
0&\zeta_2D_2(i\zeta_3\tau,x)\\
-\zeta_3D_3(i\zeta_2\tau,x)&0
\end{array}\right];\quad\widehat{D}(\tau,x)\stackrel{\rm def}{=}\left[
\begin{array}{ccc}
0&1\\
1&0
\end{array}\right]D(\tau,x);\label{eq4.25}
\end{equation}
$$\vec{A}(t,x)\stackrel{\rm def}{=}(1+a(i\zeta_2t,x),1+a(i\zeta_3t,x));\quad B(\tau,t,x)\stackrel{\rm def}{=}\left[
\begin{array}{ccc}
b_2(\tau,i\zeta_2t,x)&b_2(\tau,i\zeta_3t,x)\\
b_3(\tau,i\zeta_2t,x)&b_3(\tau,i\zeta_3t,x)
\end{array}\right].$$
Matrix $Q(t,x)$ is invertible and $\det Q(t,x)=5/4$ due to \eqref{eq3.46}. Solution of the singular integral equation is conducted with the use of a Riemann boundary value problem \cite{20, 21}.

{\bf Conclusion.} {\it From the four spectral data $\mathfrak{A}_\theta$, $\mathfrak{A}_{\widehat{\theta}}$ \eqref{eq4.2} ($\theta\not=\widehat{\theta}$) and $\mathfrak{A}_{\theta,h}$, $\mathfrak{A}_{\widehat{\theta},h}$ \eqref{eq4.12} ($h\not=0$), for which \eqref{eq4.11} holds, potential $q(x)$ is unambiguously restored.}

\begin{remark}\label{r4.1}
Scheme of solution of this inverse problem is given in Subsection 4.1.

Due to the paper volume concerns, the problem of description of spectral data $\{\mathfrak{A}_\theta,\mathfrak{A}_{\theta,h}\}$ is not considered here.
\end{remark}

\renewcommand{\refname}{ \rm 
\end{Large}
\end{document}